\documentclass[5p,times]{elsarticle} 

\usepackage[dvipsnames,table]{xcolor}

\newcommand{\cmt}[1]{{\color{black}{#1}}}
\newcommand{\cmtb}[1]{{\color{black}{#1}}}

\usepackage{amsmath}
\usepackage{siunitx}

\usepackage[ruled,linesnumbered]{algorithm2e}

\usepackage{booktabs}
\usepackage{subcaption}

\usepackage{array,makecell}

\usepackage[inline]{enumitem}

\DeclareMathOperator*{\argmax}{arg\,max}
\DeclareMathOperator*{\argmin}{arg\,min}

\usepackage{amssymb}

\journal{Additive Manufacturing}

\begin{document}

\begin{frontmatter}



\title{Stress Flow Guided Non-Planar Print Trajectory Optimization\\
for Additive Manufacturing of Anisotropic Polymers}


\author[inst1,inst2]{Xavier Guidetti}
\ead{xaguidetti@control.ee.ethz.ch}
\author[inst1]{Efe C. Balta}
\ead{ebalta@control.ee.ethz.ch}
\author[inst3]{Yannick Nagel}
\ead{yannick.nagel@nematx.com}
\author[inst4]{Hang Yin}
\ead{hanyin@student.ethz.ch}
\author[inst1,inst2]{Alisa Rupenyan\corref{cor1}}
\ead {ralisa@control.ee.ethz.ch}
\author[inst1]{John Lygeros}
\ead{lygeros@control.ee.ethz.ch}

\affiliation[inst1]{organization={Automatic Control Laboratory, ETH Zurich},
            addressline={Physikstrasse 3}, 
            postcode={8092}, 
            city={Zurich},
            country={Switzerland}}
            
\affiliation[inst2]{organization={Inspire AG},
            addressline={Technoparkstrasse 1}, 
            postcode={8005}, 
            city={Zurich},
            country={Switzerland}}
            
\affiliation[inst3]{organization={NematX AG},
            addressline={Vladimir-Prelog-Weg 5}, 
            postcode={8093}, 
            city={Zurich},
            country={Switzerland}}
            
\affiliation[inst4]{organization={Department of Mechanical and Process Engineering, ETH Zurich},
            addressline={Leonhardstrasse 21}, 
            postcode={8092}, 
            city={Zurich},
            country={Switzerland}}
            
\cortext[cor1]{Corresponding author}

\begin{abstract}
When manufacturing parts using material extrusion additive manufacturing and anisotropic polymers, the mechanical properties of a manufactured component are strongly dependent on the print trajectory orientation. We conduct non-planar slicing and optimize the print trajectories to maximize the alignment between the material deposition direction and the stress flow induced by a predefined load case. The trajectory optimization framework considers manufacturability constraints in the form of uniform layer height and line spacing. We demonstrate the method by manufacturing a load bearing mechanical bracket using a 5-axis 3D printer and a liquid crystal polymer material. The failure strength and stiffness of the optimized bracket are improved by a factor of 44 and 6 respectively when compared with conventional printing.
\end{abstract}




\end{frontmatter}


\section{Introduction} \label{sec:intro}

\begin{figure}[htbp]
     \centering
     \begin{subfigure}[t]{0.45\textwidth}
         \centering
         \includegraphics{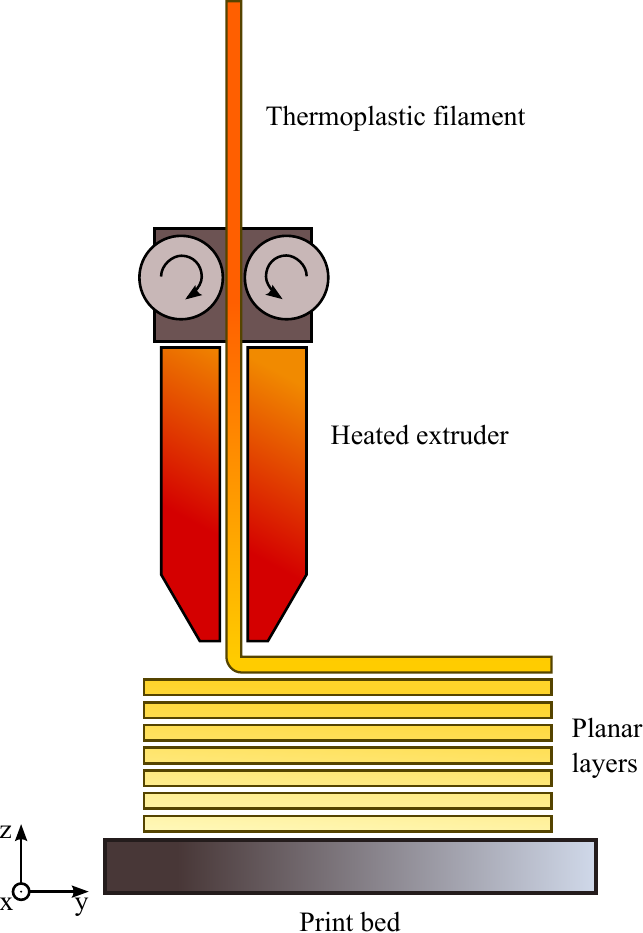}
         \caption{2.5D process}
     \label{fig:fff}
     \end{subfigure}
     \begin{subfigure}[t]{0.45\textwidth}
         \centering
         \includegraphics{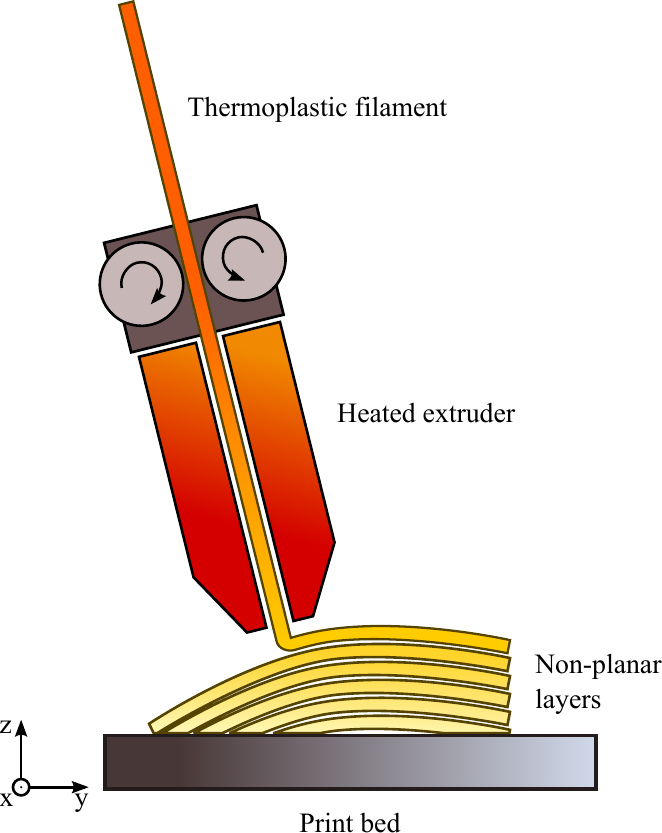}
         \caption{Non-planar process}
         \label{fig:fff-np}
     \end{subfigure}
     \caption{Simplified representation of the 2.5D and non-planar fused filament fabrication processes}
\end{figure}

In most 3D-printing applications, the material is formed in 2D layers via the printing technology of interest (e.g.\ material extrusion, laser power, or lithography). By forming subsequent layers on top of the previous ones, a 3D object is obtained. Material extrusion additive manufacturing \cite{ASTM_standard}, also known as fused filament fabrication (FFF), or fused deposition modeling (FDM), is one of the most common 3D printing processes \cite{sculpteo2021}. \cmt{Traditionally, in FFF, a thermoplastic material is deposited in predefined paths via a numerically controlled heated extruder to create a part in a layer-wise fashion.} The design geometry is sliced into planar layers that are parallel to the print bed (in the X-Y plane, see Fig.\ \ref{fig:fff}). As the final \cmt{solid} part is obtained by superposition of \cmt{planar} layers in the Z-axis, the conventional FFF approach is \cmt{often} referred to as 2.5D printing \cite{gibson2021additive}.

Based on the printing material, design geometry, printing machine, and other parameters, the design geometry is often oriented to maximize mechanical properties of interest, which may be dimensional requirements, or mechanical strength conditions~\cite{hooshmand2021optimization,di2020reliable,delfs2016optimized,wang2020simultaneous}.
While \cmt{2.5D FFF} provides material and geometry flexibility as a manufacturing process, the layer-wise processing of the parts in a fixed horizontal plane may limit the mechanical performance of the printed parts~\cite{zhang2021effect,ahn2002anisotropic}, as the inter-layer interfaces are prone to failures~\cite{gao2021fused,davis2017mechanical}.
\cmt{For example, if there is oblique loading, or if the printed part does not follow a simple planar surface, the stress flow often crosses through the inter-layer bonding locations creating premature failures. Despite the existence of methods to characterize part failure under loading \cite{capote2019failure}, the planar layer constraint in 2.5D FFF remains a limit to the applicability of the technology in mechanically demanding use cases.} 

Here instead, we study a non-planar FFF process where layers are defined by arbitrary curves and surfaces in a 3D space \cmt{(see Fig.\ \ref{fig:fff-np}). Non-planar FFF printing approaches commonly use a 5-axis printer to deposit material along non-planar trajectories, such that the stress flow on the designed geometry aligns with the beads deposition direction within layers, improving mechanical properties \cite{ramos2022optimal,fang2020reinforced,zhang2021singularity}.} While there are methods in the literature to optimize print trajectories for non-planar layers, ensuring manufacturability constraints under various geometrical features (e.g.\ holes) is often a challenge~\cite{khan2018curvilinear}. Moreover, constraining the print paths to meet conditions required by the material (e.g.\ uniform line spacing, or constant layer height) has not been addressed. The printability constraints are particularly tight for anisotropic materials. Anisotropy is however needed to exploit the potential of optimized non-planar printing~\cite{brenken2018fused}. 

\cmt{In this work, we utilize novel anisotropic polymer materials and design and manufacture non-planar print trajectories that are aligned with the stress flow induced by the loading conditions.} 
Section \ref{sec:background} discusses the existing literature and the shortcomings that are addressed by our method. In Section \ref{sec:matandmeth}, we introduce the anisotropic polymer we utilized and the details of our print optimization method. In Section \ref{sec:case_study}, the method is applied to several demonstrator and benchmarking geometries. Section \ref{sec:experiments} details the experimental results produced by manufacturing and testing one of the geometries. Finally, in Section \ref{sec:disc}, we discuss the results and compare them to the 2.5D approach.

\subsection{Background} \label{sec:background}

Multiple works have investigated the field of non-planar printing \cite{nisja2021short}. In \emph{conformal printing}, extrusion is conducted along non-planar trajectories to produce a thin shell over a non-planar substrate. This technique has been used in \cite{ahlers20193d} to finish a piece by removing the stair artifacts produced by conventional 2.5D printing. Further extending into non-planar complexity, \cite{ramos2022optimal} presents algorithms for printing complex geometries with a 5-axis printer, but the results are not demonstrated on a real system. \cmt{Similarly, print trajectory planning for non-planar robotic deposition is studied in \cite{ramos2022optimal} and in \cite{shembekar2019generating}; however both works ignore stress flow alignment and manufacturability constraints. End-to-end implementation of non-planar FFF optimization requires knowledge of the extrusion dynamics \cite{serdeczny2018experimental,aksoy2020control}, machine kinematics \cite{liniger2019real}, and efficient trajectory optimization formulations \cite{rupenyan2021performance} that consider process constraints \cite{zhang2021singularity}, material, printed geometry \cite{cicek2017numerical,balta2022numerical,balta2021layer}, and stress flow field \cite{fang2020reinforced}.}

Recent work has explored the deposition of fiber-reinforced filaments in the stress flow directions \cite{barton2019fiber,chen2022field}. Such applications are limited due to the continuity constraints on the deposited fibers and result in restricted flexibility as fibers need to be cut at the end of each individual deposition path. Other works on fiber-reinforced materials avoid the cutting problem by producing continuous trajectories \cite{yao20213d}, however without optimizing for the stress flow. \cmt{Finally, stress-orientation in 2.5D has been studied in~\cite{xia2020stress}, but only for planar stress.}

In recent research, a Liquid Crystal Polymer (LCP) material for 3D printing has been developed, patented, and adopted by NematX AG\footnote{https://nematx.com/}. \cmt{LCPs are composed of aromatic thermotropic polyesters that self-assemble into nematic domains (i.e.\ the long axes of the molecules are arranged in parallel, but not in well defined planes) when heated above their melting temperature. Alignment of different domains is achieved by extrusion through the heated nozzle of a FFF machine: the elongation and shear forces produced when the material traverses the nozzle reorient the nematic domains in the extrusion direction. The material alignment in the direction of extrusion greatly improves the mechanical properties and results in anisotropic material properties \cite{gantenbein2018three}.} By leveraging such materials, it is possible to improve the mechanical properties of the FFF printed parts, using non-planar printing methods. While this possibility has been demonstrated in customized applications, currently there exists no rigorous framework that exploits the anisotropic material properties to improve mechanical strength of FFF through non-planar printing. The quality of LCP printed parts largely depends on
\begin{enumerate*}[label=\emph{\alph*})]
    \item the quality of slicing in the desired layer thickness across all layers and how each layer supports a subsequent one (see~\cite{balta2021layer} for 2.5D printing), and
    \item the quality of trajectory generation in each layer and the homogeneous spacing of the deposited lines (see~\cite{GuidettiIFAC} for 2.5D printing).
\end{enumerate*}  
Previous approaches \cite{fang2020reinforced} have successfully tackled the optimization problem for isotropic materials such as polylactic acid (PLA), neglecting the practical constraints of high-quality FFF printing. \cmt{This resulted, for example, in toolpaths containing large variations in line spacing ($\pm 50\%$) or layer height ($\pm 30\%$).} With anisotropic polymers like LCPs that require small nozzles, low layer heights and minimal line spacing changes, these variations are too large to print a part successfully. Additionally, the application of the method from \cite{fang2020reinforced} to common load bearing brackets often resulted in print paths that are physically not printable due to their orientation.

Our novel approach aims to improve the mechanical properties of non-planar FFF printed parts by utilizing anisotropic materials such as LCPs or fiber-reinforced filaments~\cite{brenken2018fused,jiang2017anisotropic}. We formulate an optimization framework that includes the part properties (load case, material anisotropy) and the manufacturability constraints (extrusion and machine dynamics). The main contributions of this work are
\begin{itemize}
    \item A novel stress-aligned non-planar \cmt{print optimization method considering manufacturability constraints to generate trajectories that quantitatively outperform existing methods by reducing layer and line spacing deviations};
    \item A non-planar FFF printing framework suitable for mechanical design geometries with \cmt{a complex stress distribution generated by holes or elaborate mechanical features};
    \item The experimental application of stress-aligned non-planar FFF printing with anisotropic polymers to achieve $44\times$ improvement in mechanical strength over baseline approaches.
\end{itemize}

\begin{figure}[htbp]
\centerline{\includegraphics{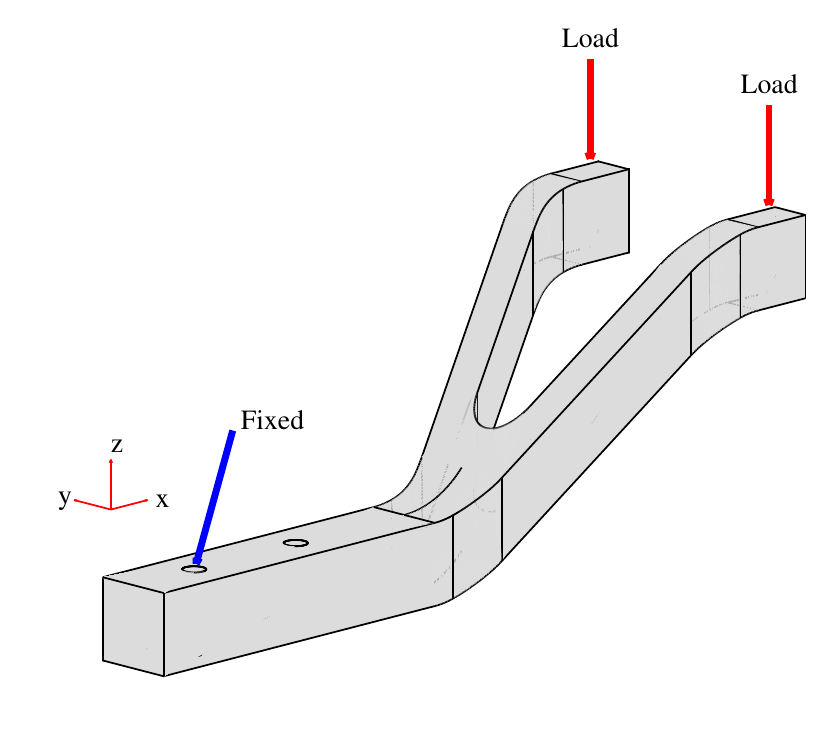}}
\caption{\emph{Fork} demonstrator geometry and loads. The part is fixed in position via screws in the holes. \cmtb{The maximum dimensions of the part are \SI{140}{\milli\metre} in the X direction, \SI{51}{\milli\metre} in the Y direction, and \SI{60}{\milli\metre} in the Z direction.}}
\label{fig:load_case}
\end{figure}

To demonstrate the slicing and print trajectory optimization methods that we propose, we will use a load bearing mechanical piece, that we refer to as the \emph{fork} for simplicity. Figure \ref{fig:load_case} shows the fork together with its load case. The loads create a non-obvious stress flow through the main body and the two arms, which makes stress-aligned filament deposition challenging. The fork is particularly challenging for the classical 2.5D approach. Irrespective of the orientation of the piece on the print bed (with or without support), slicing in planar parallel layers will not yield unbroken print lines traveling through the entire part, with an adverse impact mechanical properties. Moreover, as the main body splits in two arms, and as the different sections of the fork have different widths, it is difficult to find print trajectories leading to constantly-spaced, long, unbroken lines.

\section{Materials and methods} \label{sec:matandmeth}

\subsection{Material} \label{sec:material}

The print trajectory optimization method we propose has been developed to exploit the properties of anisotropic materials for FFF. A notable example are \cmt{LCP filaments \cite{gantenbein2018three}, developed, used for high-end applications and supplied} by NematX AG. Studies conducted on the mechanical properties of LCP parts printed with unidirectional deposition of the filaments have shown that properties are strongly dependent on the printing orientation. \cmtb{Comparison between tensile samples printed \ang{90} (i.e.\ load perpendicular to the planar print layers surfaces) and \ang{0} (i.e.\ filaments aligned to the load) to the loading direction has shown that \ang{0} printing produces a $7.0\times$ increase in the Young's modulus and a $3.6\times$ increase in the ultimate tensile strength, when compared to \ang{90} printing (see Fig.\ 3 of \cite{gantenbein2018three}).} The material is very sensitive to the parameters of the FFF deposition process \cite{GuidettiIFAC}. The best mechanical properties can be achieved when the layer height and line spacing are kept as constant as possible. Additionally, during the deposition of tightly packed lines, the drag created by the nozzle printing the next line causes misalignment in the previously deposited monomers, further reducing performance \cite{siqueira2017}.

\subsection{Methods} \label{sec:methods}

\begin{figure*}[htbp]
\centerline{\includegraphics[width=\textwidth]{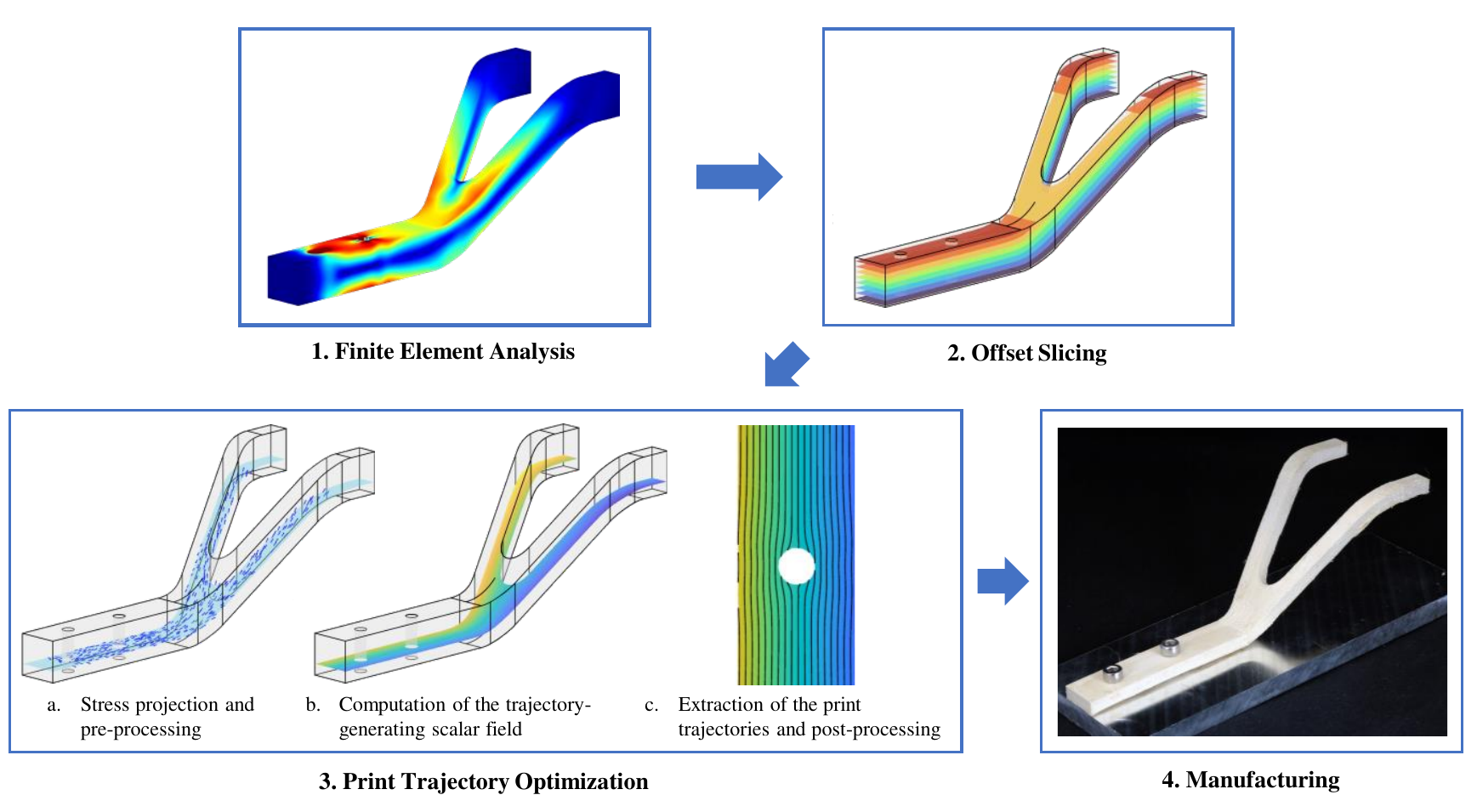}}
\caption{Graphical flowchart of the stress-aligned trajectory optimization method, demonstrated on the \emph{fork} geometry}
\label{fig:flowchart}
\end{figure*}

Figure \ref{fig:flowchart} provides a graphical overview of the proposed method for stress-aligned trajectory optimization. For a given geometry and load case, the first step is to compute the stress in the part using Finite Element Analysis (FEA) (Section\ \ref{sec:fea}). We then divide the geometry in constantly spaced non-planar slices. The orientation of the slices is selected to reduce inter-layer stress, allow long and unbroken fibers across the entire geometry, and maximize print quality through uniform layer height (Section\ \ref{sec:slicing}). The third step is to fill each non-planar layer with optimized print trajectories. We first project the stress computed in FEA on the slice surface, and pre-process it so that it matches the structure required by the optimization problem. Solving the problem corresponds to computing a trajectory-generating scalar field in each layer. The isostatic lines (isolines) of this field -- which we post-process for manufacturability -- are the stress-aligned and homogeneously spaced print trajectories that form the layer (Section\ \ref{sec:traj_opt}). Once each layer has been optimized, we join all the trajectories in an ordered sequence of points that is sent to a 5-axis FFF machine to manufacture the part. In the next sections, we explain each step of the method in detail.

\subsubsection{Finite Element Analysis} \label{sec:fea}

To generate stress-aligned print trajectories, we must first compute the stress in the parts we intend to manufacture. This is achieved using FEA to simulate the specific load case for which a piece is going to be optimized \cite{reddy2019FEA}. We begin by subdividing the 3D representation of a part in a tetrahedral mesh. After having applied the boundary conditions corresponding to the forces on the part and having assigned to the part isotropic properties (i.e.\ an arbitrary Young's modulus and a Poisson's ratio $\nu = 0.3$), we solve the resulting discretized problem leading to the Cauchy stress tensor $\boldsymbol\sigma$ for each node in the mesh. A general stress tensor is typically decomposed using eigenvalue decomposition to obtain the principal stresses and principal directions \cite{bower2009applied}. In a reference frame oriented along the principal directions, the stress tensor can be represented as a diagonal matrix where the diagonal entries $\sigma_{1}$, $\sigma_{2}$ and $\sigma_{3}$ are ordered so that \cmt{$|\sigma_{1}| \geq |\sigma_{2}| \geq |\sigma_{3}|$. We call $\sigma_{1}$ the \emph{maximum principal stress} and $\sigma_{3}$ the \emph{minimum principal stress} and the corresponding normalized eigenvectors \emph{maximum principal stress direction} and \emph{minimum principal stress direction}.} They represent the direction in which the principal stresses act. The maximum and minimum principal stress directions are orthogonal. We will call \emph{principal stress vectors} the vectors obtained by multiplying a principal stress with the corresponding principal stress direction and \emph{stress flow} the vector field that regroups the principal stress vectors computed at all the nodes.

\subsubsection{Offset Slicing} \label{sec:slicing}

Conventionally, in 2.5D printing the desired component is sliced in evenly spaced layers that are parallel to the print bed. The homogeneity of the layers allows one to obtain high-quality prints with excellent layer adhesion, producing mechanically strong pieces. We adapt this slicing scheme to the non-planar context, where high-degrees-of-freedom machines can extrude along curved layers. As in 2.5D printing, we choose as our first slice the contact surface between the part and its (curved) support. Then, we generate the following slices so that each one has a constant distance from the one below. Given a desired geometry and having chosen which one of its faces will be resting on a support piece, we select this entire face as our first slice. Then, we follow the geodesic heat method (Alg. \ref{alg:heat_method}) proposed in \cite{Crane_Geodesics} to compute the geodesic distance of a domain (the desired geometry) from a subset (the first slice).

\begin{algorithm}[htbp]
\caption{The Geodesic Heat Method \cite{Crane_Geodesics}}\label{alg:heat_method}
Compute the flow of heat $u$ from a contact surface, by integrating the heat flow $\dot{u} = \Delta u$ for a fixed small time $t$; \\
From the temperature gradient $\nabla u$, compute the normalized and negated vector field $X = -\nabla u/|\nabla u|$; \\
Recover the final distance $\phi$, whose gradient follows X, by solving the Poisson equation $\Delta\phi = \nabla \cdot X$. \\
\end{algorithm}

The practitioner must first select the orientation of the piece on the print bed or support. Clearly, when working with the objective to maximize the alignment between the deposited fibers and the maximum principal stress flow, we would like the normals to the print layers to be aligned with the minimum principal stress flow. Existing methods \cite{fang2020reinforced} that generate slices maximizing the alignment between layer normals and minimum principal stress vectors simply push the orientation problem downstream, forcing the practitioner to look for ways to manufacture irregularly oriented layers whose complexity leads to poor prints. Often, this results in twisted or intertwined layers which cannot be physically manufactured. As our method focuses on manufacturability, the initial choice of orientation is actually an advantage, as it allows one to select the orientation producing the best print quality. Additionally, when printing mechanical pieces designed for a well defined load case, we noticed it is fairly easy and intuitive to find an initial orientation generating a satisfactory layer alignment. We support this claim empirically in Section \ref{sec:impl_slicing}.

\subsubsection{Print Trajectory Optimization} \label{sec:traj_opt}

The starting point of our approach to print trajectory optimization is the work of Fang et al. \cite{fang2020reinforced}, that we embellish with numerous refinement steps which we discuss in detail. We first note that irrespective of the slicing method, there will be regions where the maximum principal stress vector is not tangent to the surface of the slice. Thus, we begin by projecting the maximum principal stress flow field onto the slice. Because of the structure of the print path optimization problem, we then need to rotate the projected vectors around each surface normal to obtain their orthogonal vectors $\mathbf{f_\perp}$. For a vector $\mathbf{f}$ (belonging to the stress flow $\mathbf{F}$) and the corresponding surface normal $\mathbf{n}$, we compute $\mathbf{f_\perp}$ as
\begin{equation}
    \mathbf{u} = \mathbf{f} - \frac{\mathbf{f}\cdot\mathbf{n}}{\mathbf{n}\cdot\mathbf{n}}\mathbf{n}\,, \quad \mathbf{f_\perp} = \frac{\mathbf{u} \times \mathbf{n}}{\|\mathbf{u} \times \mathbf{n}\|}\,.
\end{equation}
We regroup all calculated $\mathbf{f_\perp}$ in the projected, orthogonal and normalized maximum principal stress flow $\mathbf{F_\perp}$.

\begin{figure*}[htbp]
     \centering
     \begin{subfigure}[b]{0.23\textwidth}
         \centering
         \includegraphics[scale=1]{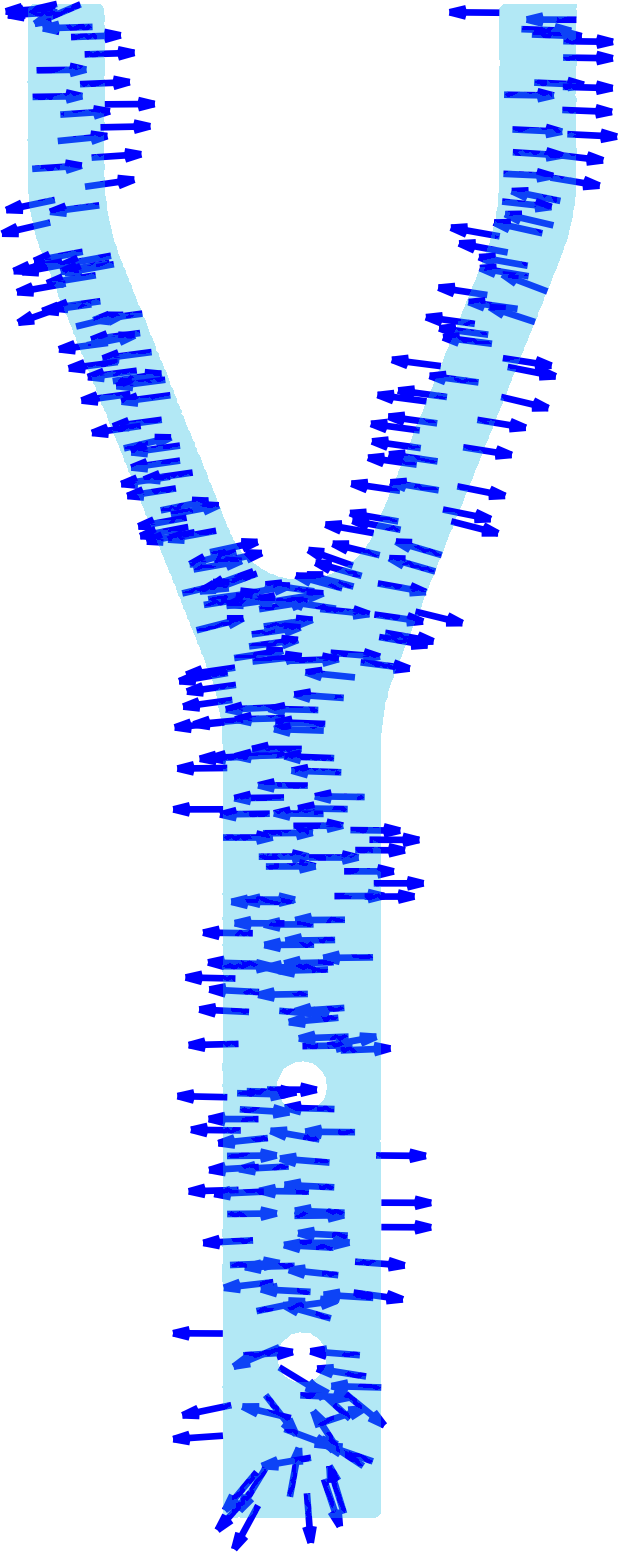}
         \caption{Before stress flow rectification}
     \end{subfigure}
     \begin{subfigure}[b]{0.23\textwidth}
         \centering
         \includegraphics[scale=1]{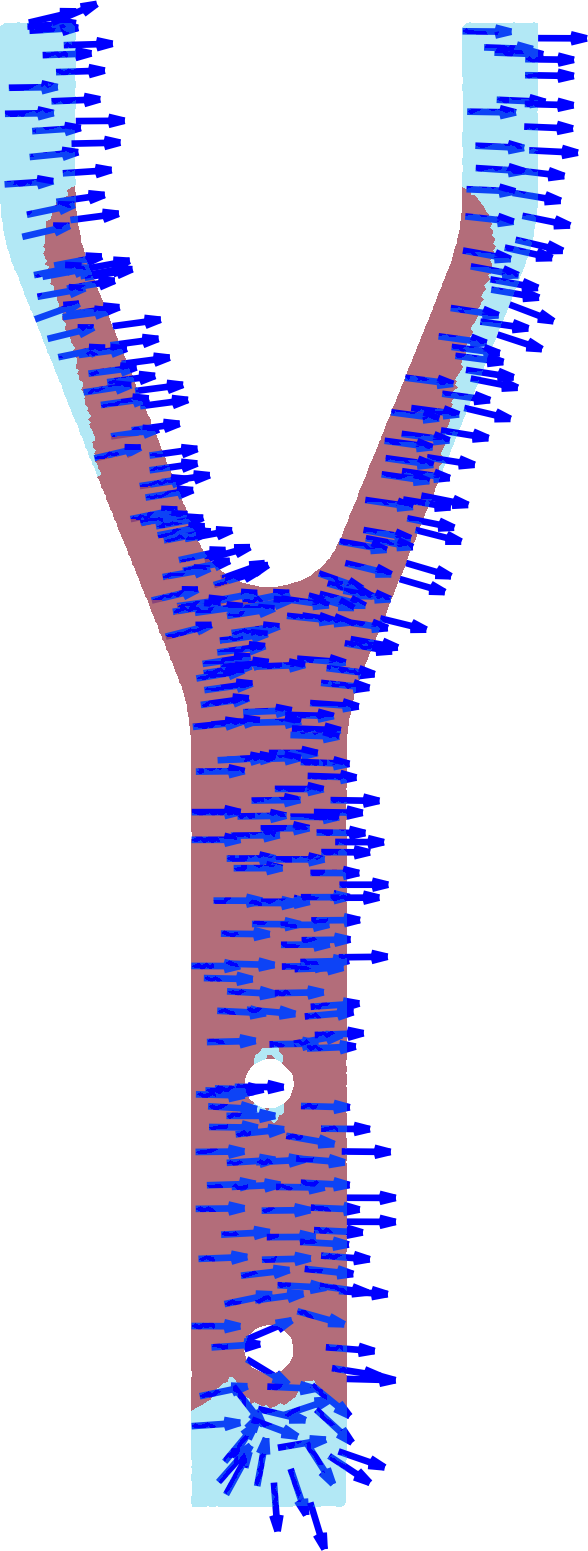}
         \caption{After stress flow rectification}
     \end{subfigure}
     \begin{subfigure}[b]{0.23\textwidth}
         \centering
         \includegraphics[scale=1]{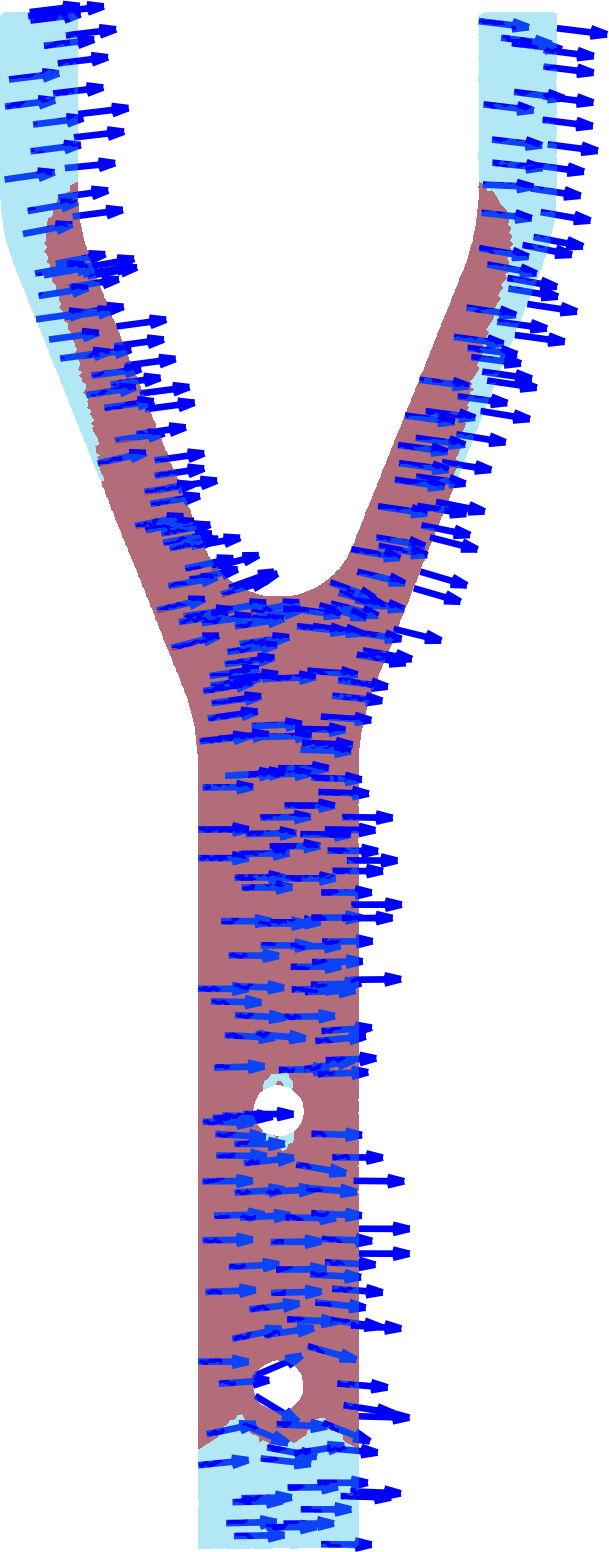}
         \caption{After stress flow propagation}
     \end{subfigure}
        \caption{Rectification and propagation of the stress flow, shown on a slice of the \emph{fork} demonstrator geometry for the load case that will be discussed in Sec.\ \ref{sec:case_study}. Note the agreement in the vectors orientation after rectification and the change of stress flow direction in the uncritical region after rectification.}
        \label{fig:stress_prop}
        \label{fig:stress_rect}
\end{figure*}

The heterogeneous stress orientations produced by the FEA, coupled with the fact that surface normals can lie on both side of the slice, can create ambiguities in the stress flow $\mathbf{F_\perp}$. This undesired effect occurs in regions where the vectors composing the stress flow are aligned (approximately) in the same direction, but have opposite orientations. As the toolpath optimization problem that we pose is sensitive to vector orientation ambiguity, we propose to pre-process the stress flow. We first analyze all the vectors in $\mathbf{F_\perp}$ to determine their \emph{main} Cartesian component. We compute
\begin{equation}
    i_d = \argmax_i{\left\|\mathbf{F_\perp}\cdot \mathbf{e}_i\right\|}\,,
\end{equation}
where $i_d$ represents the Cartesian component to which the vectors in $\mathbf{F_\perp}$ have the best alignment and $\textstyle{\mathbf{e}_i,~i=1,2,3}$ are the unit axes of a Cartesian reference frame in $\mathbb{R}^3$. Now, we can rectify the stress flow by calculating
\begin{equation}
    \mathbf{f_r} =
    \begin{cases}
      \mathbf{f_\perp} & \text{if}\ \mathbf{f_\perp}\cdot \mathbf{e}_{i_d} \geq 0 \\
      -\mathbf{f_\perp} & \text{otherwise}
    \end{cases}
\end{equation}
for all vectors in $\mathbf{F_\perp}$. We then regroup the resulting $\mathbf{f_r}$ in the rectified stress flow $\mathbf{F_r}$. Figure \ref{fig:stress_rect} includes an example of stress flow rectification. Intuitively, the rectification consists in flipping the vectors that do not lie in main stress flow orientation. The rectification ensures that the print trajectories that we then generate are homogeneous and easily manufacturable. We point out that our method requires the stress flow to be rectified, but is not sensitive to the rectified vectors orientation. This is because the extruded material should simply align with the main traction/compression direction, making the problem invariant to $\ang{180}$ rotations of the stress flow. 

\begin{figure}[htbp]
     \centering
     \begin{subfigure}[b]{0.23\textwidth}
         \centering
         \includegraphics[scale=1]{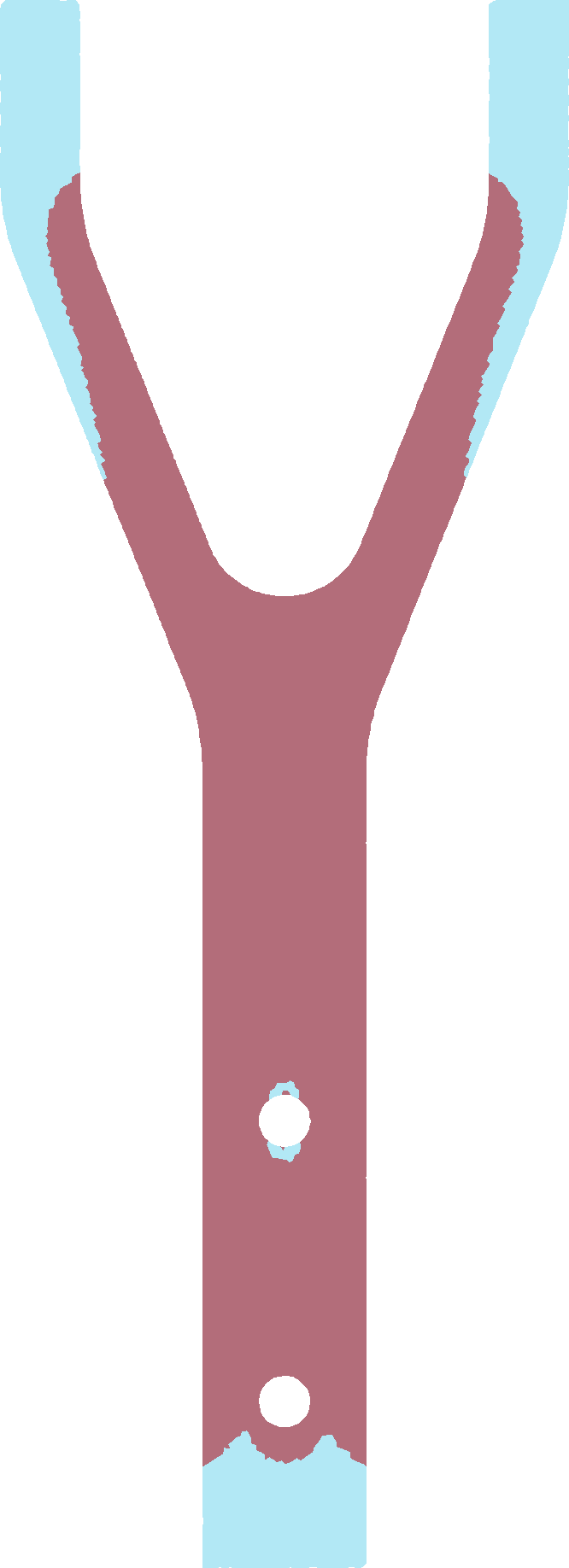}
         \caption{$\theta_a = 3 \, ,\, \theta_s = 0.1$}
         \label{subfig:crit_sel}
     \end{subfigure}
     \hfill
     \begin{subfigure}[b]{0.23\textwidth}
         \centering
         \includegraphics[scale=1]{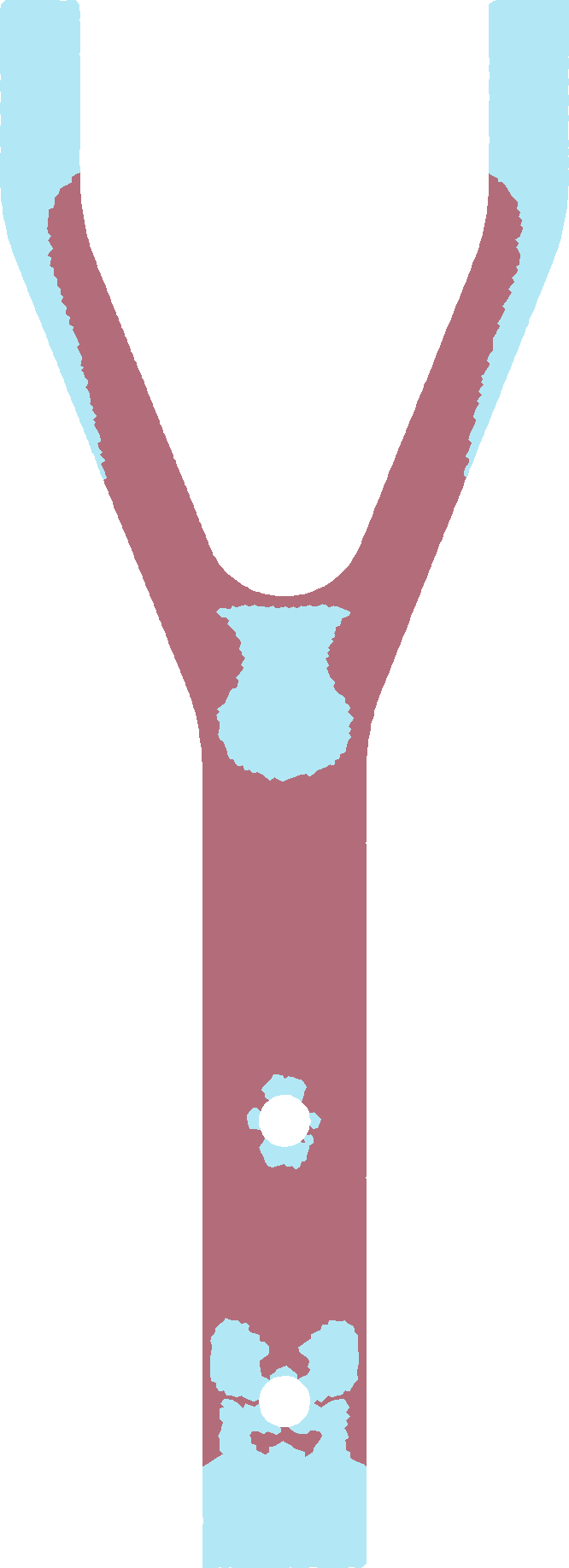}
         \caption{$\theta_a = 15 \, ,\, \theta_s = 0.1$}
     \end{subfigure}
     \hfill
     \vspace{5pt}
     \begin{subfigure}[b]{0.23\textwidth}
         \centering
         \includegraphics[scale=1]{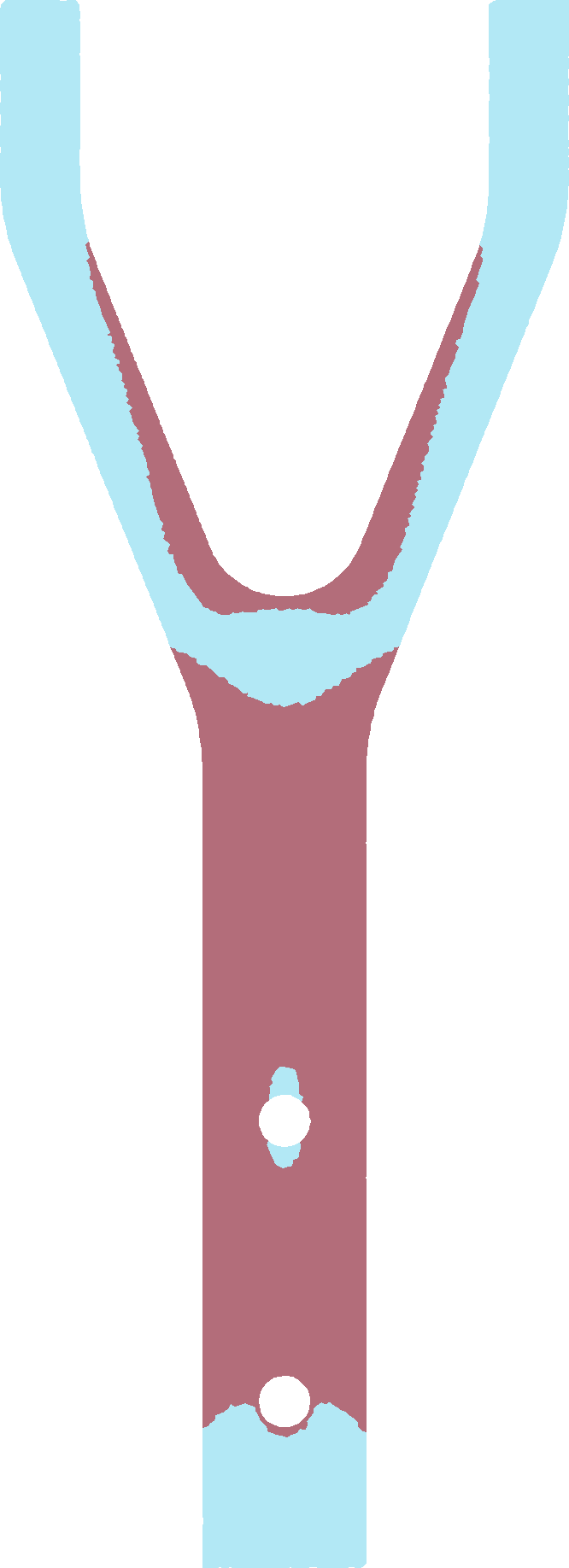}
         \caption{$\theta_a = 3 \, ,\, \theta_s = 0.2$}
     \end{subfigure}
     \hfill
     \begin{subfigure}[b]{0.23\textwidth}
         \centering
         \includegraphics[scale=1]{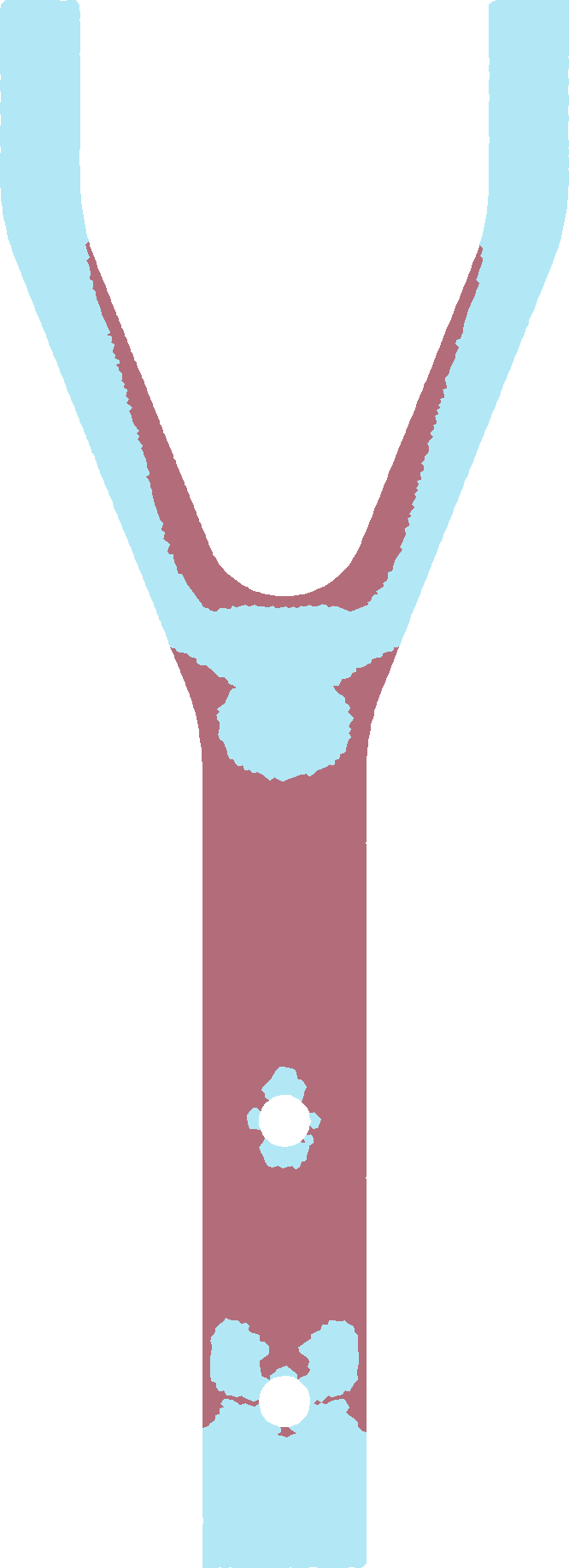}
         \caption{$\theta_a = 15 \, ,\, \theta_s = 0.2$}
     \end{subfigure}
        \caption{Critical region shape (in red) as a function of the stress anisotropy and significance parameters $\theta_a$ and $\theta_s$, shown on a slice of the \emph{fork} demonstrator geometry for the load case that will be discussed in Sec.\ \ref{sec:case_study}. The same example slice will be used in all other figures. Experimentation suggests that the method produces better trajectories when selecting parameters that produce a \emph{simply connected} critical region, which in this scenario corresponds to case (\subref{subfig:crit_sel}).}
        \label{fig:crit_comparison}
\end{figure}

Next we distinguish \emph{critical} stress regions from \emph{uncritical} ones. We consider a node as critical if the stress is both sufficiently anisotropic and sufficiently large. The former condition translates to 
\begin{equation}
\cmt{
    \frac{|\sigma_1|}{|\sigma_3|} > \theta_a\,,}
    \label{eq:anis_crit}
\end{equation}
where $\theta_a$ is the anisotropy threshold; the latter to 
\begin{equation}
\cmt{
    \frac{|\sigma_1|}{\max_\text{nodes}|\sigma_1|} > \theta_s\,,}
    \label{eq:stress_crit}
\end{equation}
where $\theta_s$ is the stress significance threshold.
The idea behind this classification is that some regions of the desired geometry undergo large and quasi unidirectional stresses: in these regions the alignment between deposited material and the stress flow is crucial. Other regions are subject to low or isotropic stresses and the material alignment is not important. In uncritical regions, we can print in the most convenient direction, which is typically the continuation of the print direction used in the critical regions. This allows us to deposit long unbroken fibers in our manufactured piece, fully exploiting LCPs properties. Figure \ref{fig:crit_comparison} illustrates on an example geometry the critical regions corresponding to different choices of the stress anisotropy and significance parameters.

It is now possible to extrapolate the stress flow $\mathbf{f_r}$ of the nodes at the critical region boundary to the neighboring uncritical nodes. We propagate the critical stress flow, until all the uncritical nodes have their original stress vector replaced by an extrapolated (and normalized) critical stress vector. The extrapolation is conducted by solving an optimization problem that minimizes the Dirichlet energy of the stress flow in the uncritical region, constrained by the boundary conditions with the critical region. The optimal extrapolated stress flow is the minimizer of the constrained problem
\begin{equation}\label{eq:dirichlet}
    \begin{array}{rll}
        \mathbf{\hat{F}_{uc}} \, = \, \argmin_\mathbf{F_{uc}} & \| G \mathbf{F_{uc}} \| ^2 \\
        \text{s.t.} & \mathbf{f_{uc}} = \mathbf{f_{r}} &  \text{for each node at the boundary} \\
        & \mathbf{f_{uc}} \cdot \mathbf{n} = 0 & \text{for each node in the uncritical region} \,,
    \end{array}
\end{equation}
where $\mathbf{F_{uc}}$ regroups all the stress flow vectors $\mathbf{f_{uc}}$ of the uncritical regions and $G$ is the numerical gradient over the triangular mesh of the uncritical region. Solving \eqref{eq:dirichlet} corresponds to finding the nodal vectors that minimize the orientation variability in the entire vector field. The two constraints enforce the requirement that the vectors at the critical-uncritical boundary are fixed and that the solutions found in the uncritical region are tangent to the surface. The solution $\mathbf{\hat{F}_{uc}}$ replaces the stress flow of the uncritical nodes. Figure  \ref{fig:stress_prop} includes an example of the propagation process.

Having applied all the previous steps to the results of the FEA simulation, to each node in the slice corresponds a pre-processed orthogonal stress vector $\mathbf{f_p}$. We regroup all $\mathbf{f_p}$ in $\mathbf{F_p}$. We are ultimately interested in solving the problem 

\begin{equation} \label{eq:main_opt}
    \hat{\phi} = \argmin_\phi \int_S \| \nabla \phi - \mathbf{F_p} \| ^2 ds \, , 
\end{equation}

to find a scalar field $\hat{\phi}$ whose gradient matches as closely as possible the orthogonal stress flow $\mathbf{F_p}$ on the slice $S$. Given the discretized nature of our problem, we reformulate the optimization problem \eqref{eq:main_opt} as the following regularized optimization problem
\begin{equation} \label{eq:least_sq_reg}
    \hat{\phi} = \argmin_\phi \| G \phi - \mathbf{F_p} \| ^2 + \epsilon \| \phi \|^2,
\end{equation}
where $G$ is the numerical gradient over the triangular mesh of the slice and $\epsilon$ is the coefficient of the regularizer.
There exists a closed form solution of \eqref{eq:least_sq_reg},
\begin{equation}
    \hat{\phi} = (G^\top G + \epsilon I)^{-1}G^\top \mathbf{F_p} \, .
    \label{eq:scalar_field}
\end{equation}
We can now produce the optimized print trajectories, which are isolines of the scalar field $\hat{\phi}$ that we just obtained. Calculating the isolines is a straightforward operation that just requires extracting equally spaced isolines on the slice mesh, based on $\hat{\phi}$ . The isolines spacing will correspond to the distance between two print trajectories. We conduct a final post-processing on the print trajectories by cubic spline re-sampling \cite{wolberg1988cubic}. We first fit a cubic smoothing spline interpolation, which produces a quasi-identical trajectory, but reduces major roughness. The spline is then sampled at even distances to produce a sequence of points that represent the print trajectory. This final smoothing step improves manufacturability and ensures that the mechanical properties of the part will not be limited by a poor quality print.

\section{Calculation} \label{sec:case_study}

In the following sections, we apply the slicing and trajectory optimization methods that we have just described to the \emph{fork} demonstrator geometry. The entire implementation is conducted using MATLAB. We first mesh the part, apply the loads shown in Fig.\ \ref{fig:load_case} and conduct FEA to compute the Cauchy stress tensors that will be used in slicing and trajectory optimization.

\subsection{Offset Slicing} \label{sec:impl_slicing}

\begin{figure*}[htbp]
    \begin{subfigure}[t]{0.49\textwidth}
        \centering
        \includegraphics{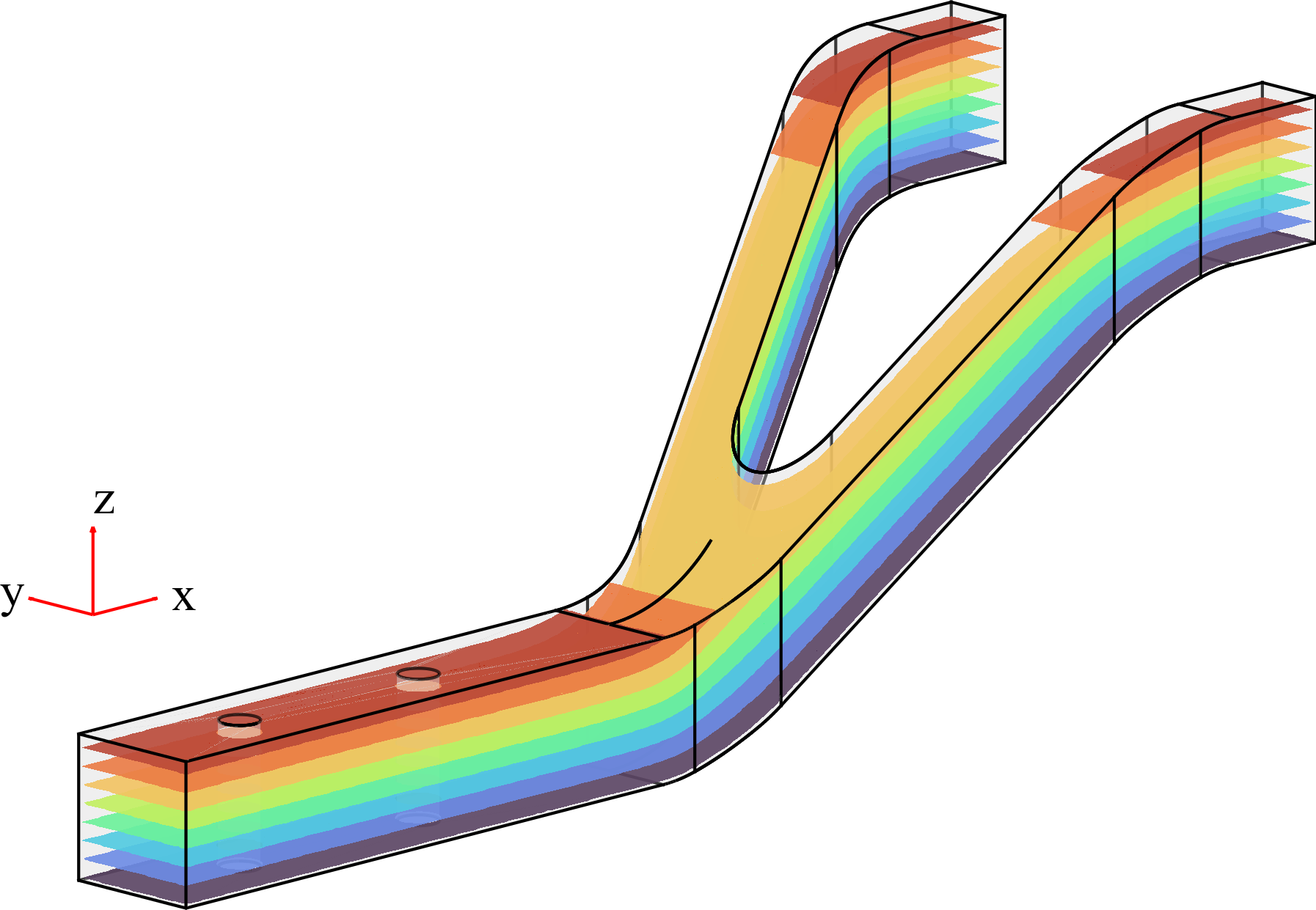}
        \caption{Offset slicing produced by our method}
        \label{fig:slice_demo_offset}
    \end{subfigure}
    \hfill
    \begin{subfigure}[t]{0.49\textwidth}
        \centering
        \includegraphics{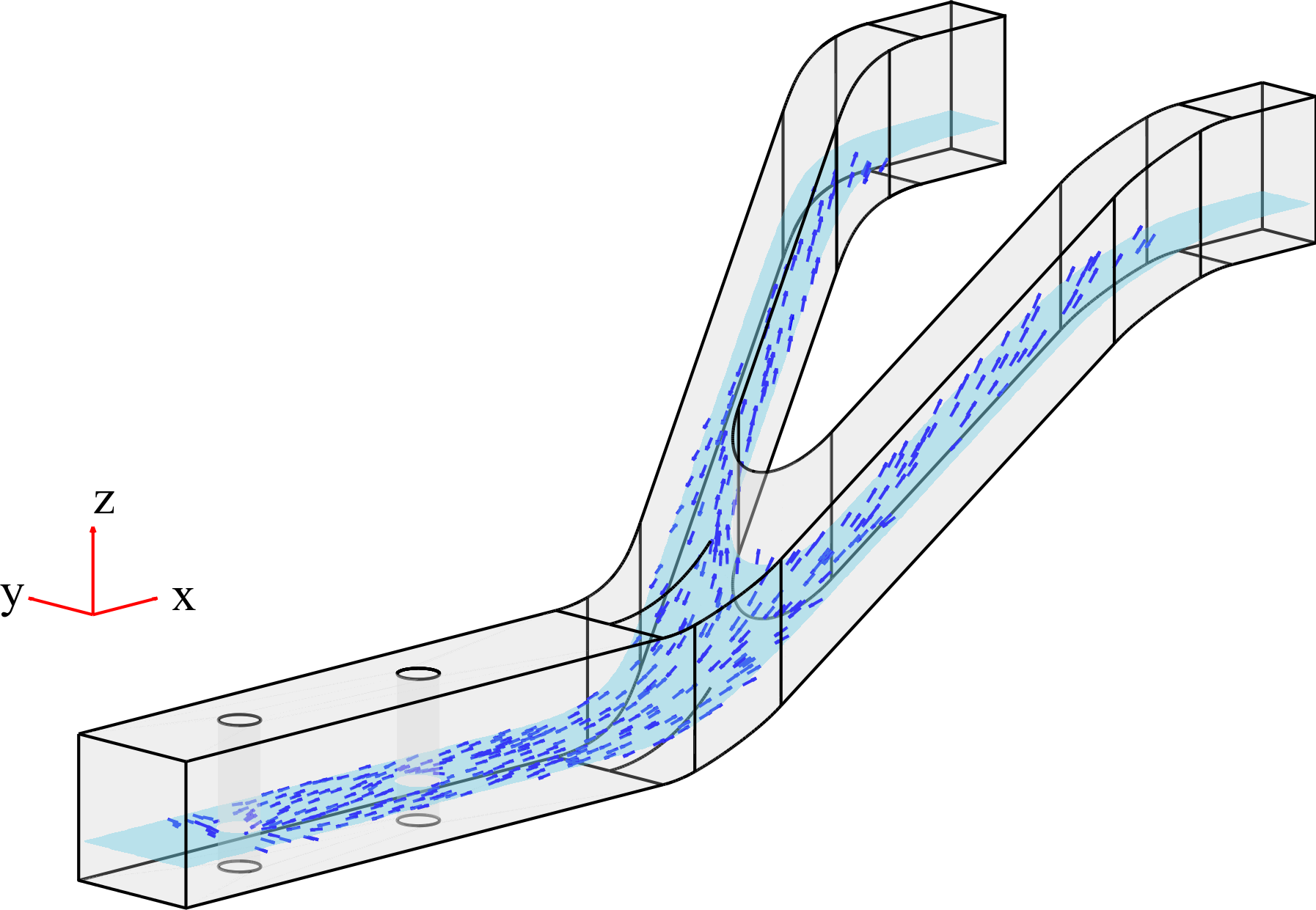}
        \caption{Example of one slice obtained via offset slicing (light blue) with the critical maximum principal stress directions (dark blue) belonging to the slice. Note how most of the vectors are approximately tangent to the slice surface, producing an average slicing stress-alignment $\bar{\gamma} \approx 1$.}
        \label{fig:slice_stress}
    \end{subfigure}
 
    \begin{subfigure}[t]{0.49\textwidth}
        \centering
        \includegraphics{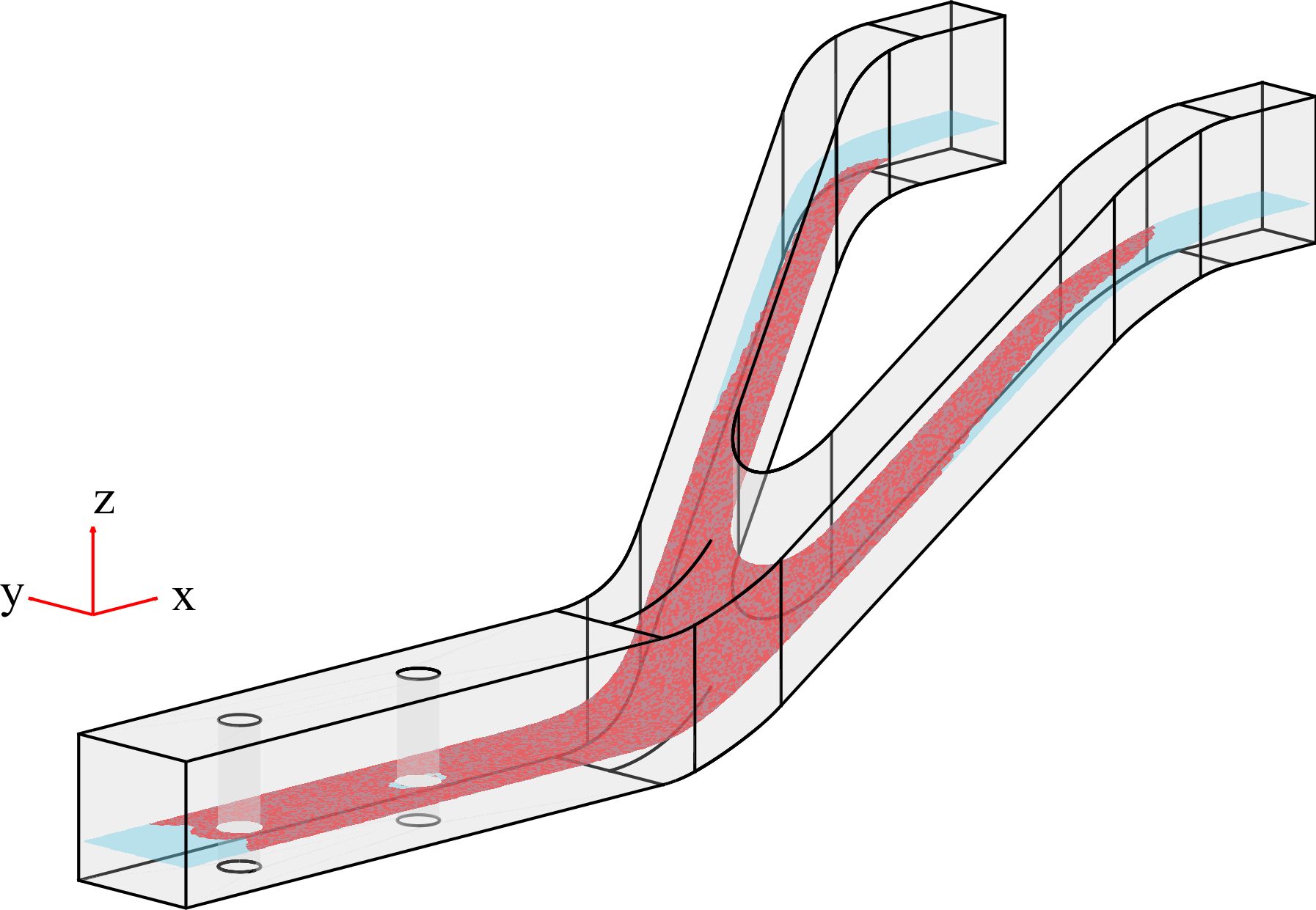}
        \caption{Critical region (in red) on one slice of the part. The critical region forms a \emph{simply connected} surface and excludes the areas where the stress is low or isotropic.}
        \label{fig:crit_region}
    \end{subfigure}
    \hfill
    \begin{subfigure}[t]{0.49\textwidth}
        \centering
        \includegraphics{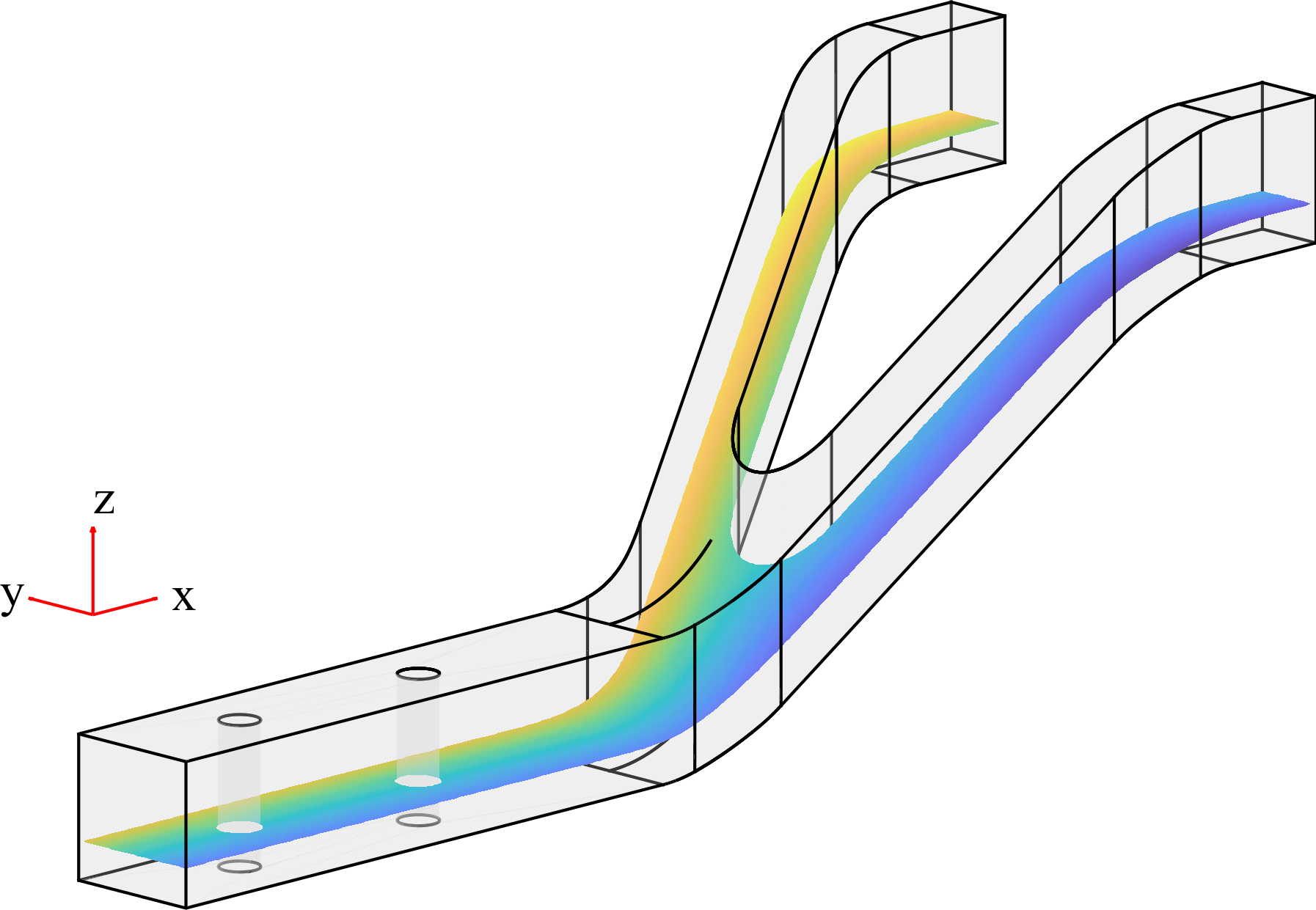}
        \caption{Trajectory-generating scalar field on one slice of the part. The color map is associated to the value of the scalars associated to each node, showing a smooth gradient.}
        \label{fig:scalar_field}
    \end{subfigure}
    
    \begin{subfigure}[t]{\textwidth}
        \centering
        \includegraphics[scale=0.7,angle=-90]{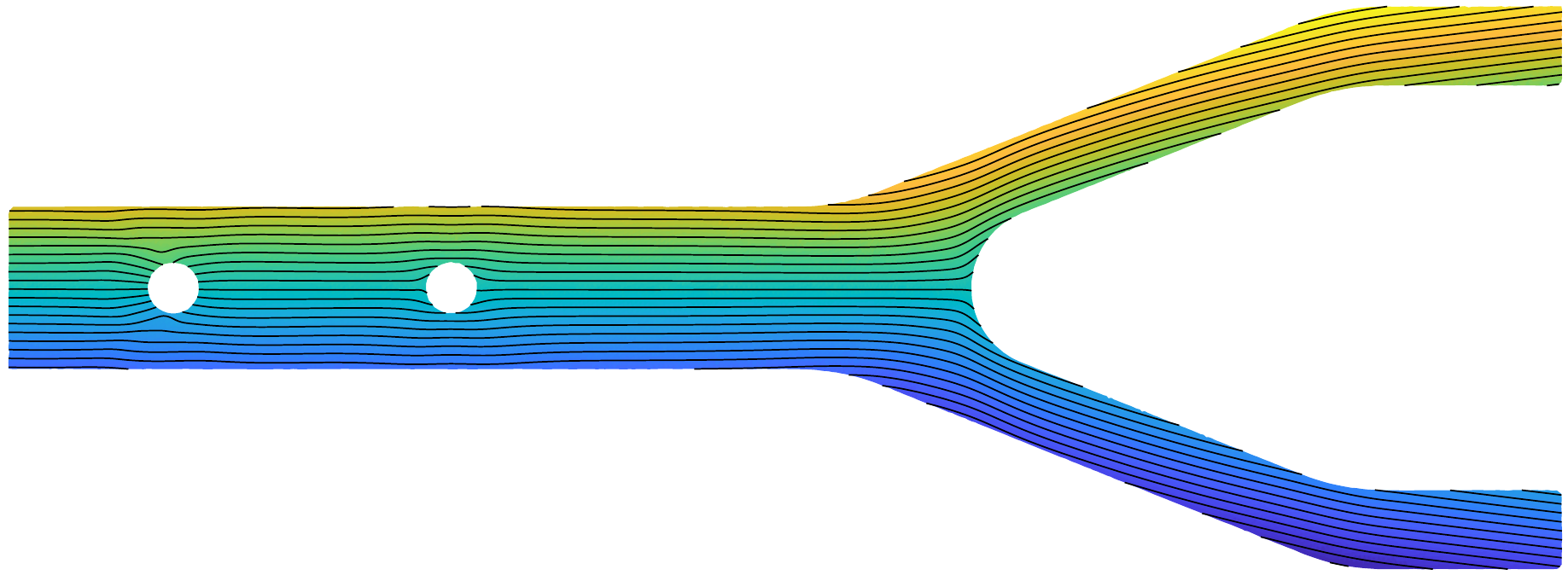}
        \caption{Isolines extracted from the trajectory-generating scalar field on one slice of the part. Note how the lines are very uniformly distributed, ensuring that the trajectories can produce a high quality print. Additionally, the trajectories follow the arms of the bracket closely and connect the extremities of the parts with unbroken print lines.}
        \label{fig:isolines}
    \end{subfigure}
    \caption{Main steps of the implementation on the \emph{fork} geometry}
\end{figure*}

For each node of the FEA mesh, we compare the the maximum principal stress vectors with the slicing surfaces normals. Following Sec.\ \ref{sec:traj_opt}, the ideal slicing is obtained when each maximum principal stress is perpendicular to the corresponding slice normal. Given the maximum principal stress vector $\mathbf{f}$ and the slicing surface normal $n$ (both assumed to be normalized) at the node $j$, we compute the nodal \cmt{slicing} stress-alignment as
\begin{equation}
    \gamma_j = \|\mathbf{f} \times \mathbf{n}\| \,.
\end{equation}
We calculate the average slicing stress-alignment over the entire piece by considering the \emph{critical} regions only. For a meshed piece containing $J$ critical nodes, we obtain
\begin{equation}
    \bar{\gamma} = \frac{1}{J} \sum_{\text{critical } j} \gamma_j\,.
    \label{eq:gamma}
\end{equation}
The proposed metric will produce $\bar{\gamma} \in [0,1]$, where $\bar{\gamma}=1$ corresponds to a perfectly aligned slicing and $\bar{\gamma}=0$ to ta case where the maximum principal stress vector is perpendicular to the slices. A comparison of our method (Fig.\ \ref{fig:slice_demo_offset}) with standard planar slicing along the X, Y, and Z axes of the demonstrator bracket reference frame is quantified in Table \ref{tab:slice_alignment}.
\begin{table}[htbp]
\centering
\caption{Average slicing stress-alignment for different slicing methods (evaluated in the critical regions only), with identical layers distance}
\renewcommand{\arraystretch}{1.3}
\begin{tabular}[h]{@{}l r r r r @{}}
\toprule
& Offset & Planar X & Planar Y & Planar Z\\
\midrule
$\bar{\gamma}$ & 0.96 & 0.38 & 0.98 & 0.87 \\
\bottomrule
\end{tabular}
\label{tab:slice_alignment}
\end{table}
Both the offset slicing method we propose and the planar slicing along the Y-axis have very high alignment scores $\bar{\gamma}$. The Y-axis planar slicing, however, produces a large fraction of slices that belong exclusively to the fork arms. As a consequence, the manufactured part would be more prone to inter-layer fractures in the arms. Conversely, the proposed offset slicing allows to work on layers that span the entire piece, significantly increasing the strength of the printed part and making offset slicing superior to all planar slicing methods.
Figure \ref{fig:slice_stress} shows one of the slices obtained using offset slicing, together with the critical maximum principal stress directions that were used to compute $\gamma$ on the depicted slice.

\subsection{Print Trajectory Optimization}

\begin{figure}[htbp]
     \centering
     \begin{subfigure}[b]{0.49\textwidth}
         \centering
         \includegraphics{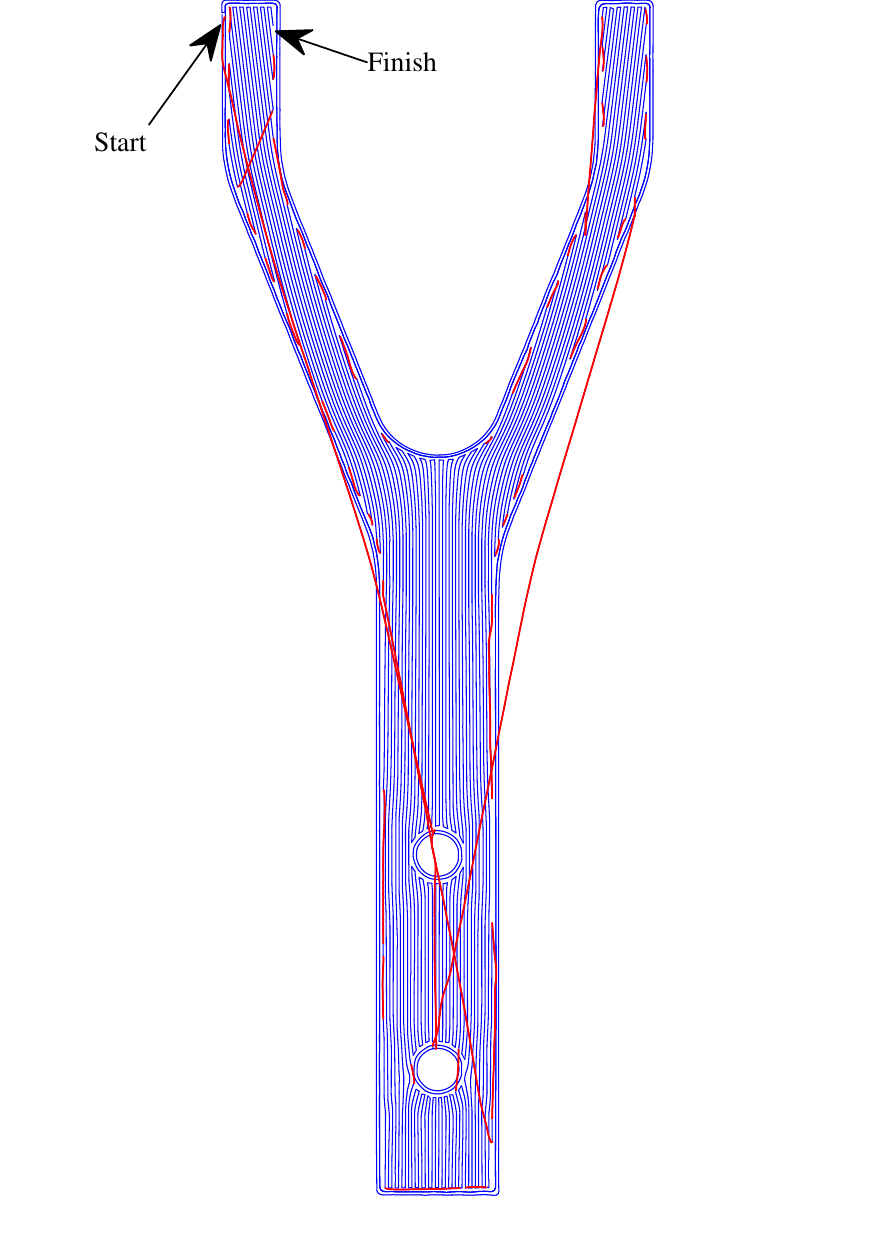}
     \end{subfigure}
     \vskip\baselineskip
     \begin{subfigure}[b]{0.49\textwidth}
         \centering
         \includegraphics{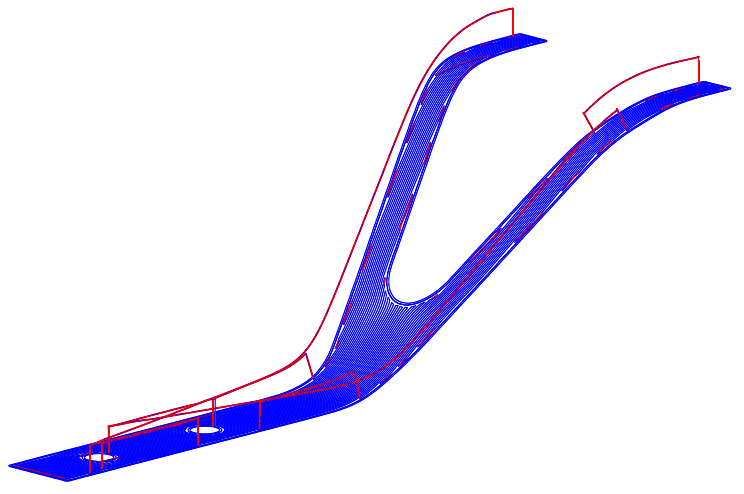}
     \end{subfigure}
        \caption{Post-processed print trajectories on one layer of the part. Travel moves are depicted in red.}
        \label{fig:final_traj}
\end{figure}

Next, we generate the print trajectories for each slice individually. We first project and rectify the stress flow on a slice, and then determine the critical regions. We set the parameters of Equations \eqref{eq:anis_crit} and \eqref{eq:stress_crit} to $\theta_a = 3$ and $\theta_s = 0.1$, to produce a \emph{simply connected} critical region that covers most of each slice. The following step consists in solving Eq. \eqref{eq:scalar_field} to compute the values of the trajectory-generating scalar field over the nodes of each meshed slice. We show the critical region and the scalar field corresponding to the example slice in Fig.\ \ref{fig:crit_region} and \ref{fig:scalar_field}.
Finally, it is possible to extract the isolines of the scalar field that was just calculated. The isolines spacing corresponds to the distance between two neighboring print trajectories, and can be easily selected to match the material and printer requirements. Figure \ref{fig:isolines} shows the isolines that were obtained from the scalar field of Fig.\ \ref{fig:scalar_field}.
In order to maximize the quality of the print, we post-process the isolines. These are first interpolated using a cubic spline with a smoothing parameter $p = 0.95$, to smoothen the small-scale wiggles that can be caused by meshing and numerical approximation. Then we add contour lines to each slice and trim the smoothed isolines to avoid overlap with the contours during printing. We finally connect sequentially the adjacent print lines and generate travel moves where required. Figure \ref{fig:final_traj} shows the final print trajectory that will be sent to the printer to manufacture one layer.

\begin{figure*}[htbp]
     \centering
     \begin{subfigure}[t]{0.49\textwidth}
         \centering
         \includegraphics{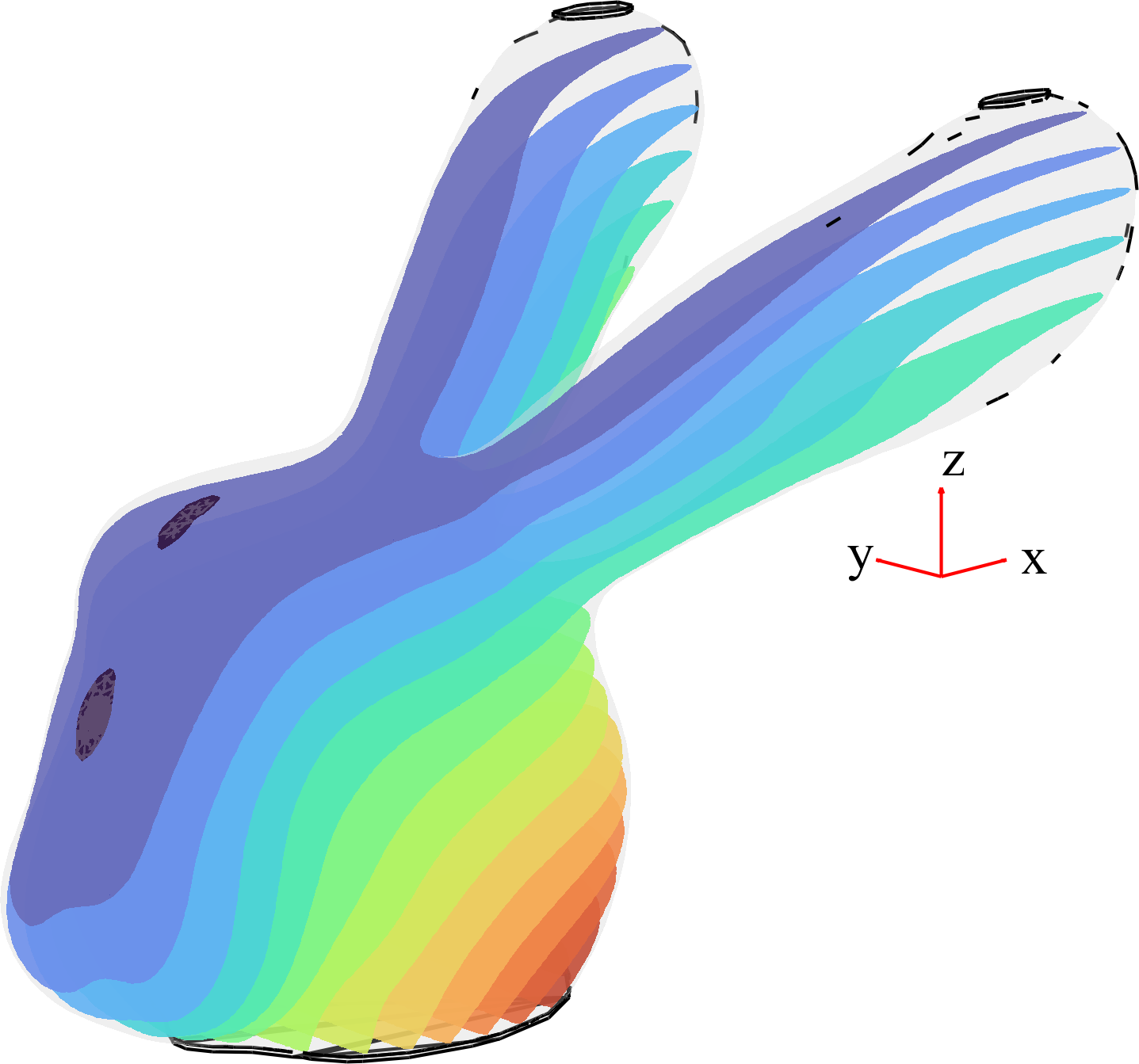}
         \caption{Offset slicing}
     \end{subfigure}
     \hfill
     \begin{subfigure}[t]{0.49\textwidth}
         \centering
         \includegraphics{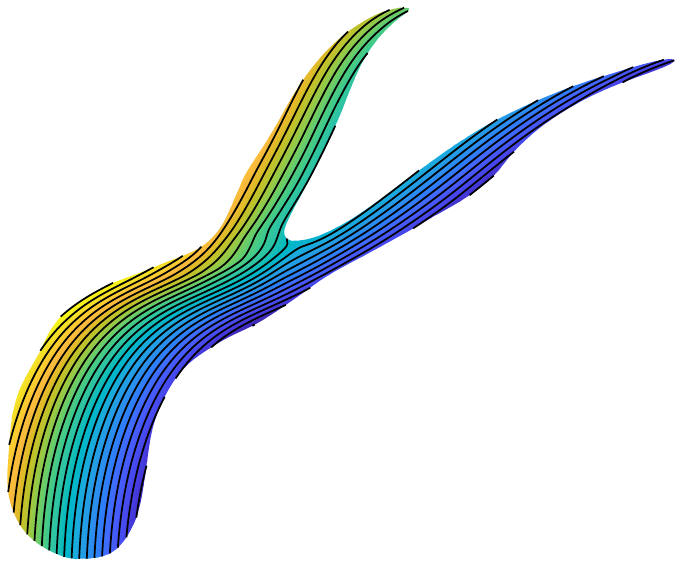}
         \caption{Isolines and scalar field on one layer of the part}
     \end{subfigure}
        \caption{Application of our method to the benchmarking piece known as \emph{bunny head}}
        \label{fig:bunny}
\end{figure*}

\begin{figure*}[htbp]
     \centering
     \begin{subfigure}[t]{0.49\textwidth}
         \centering
         \includegraphics{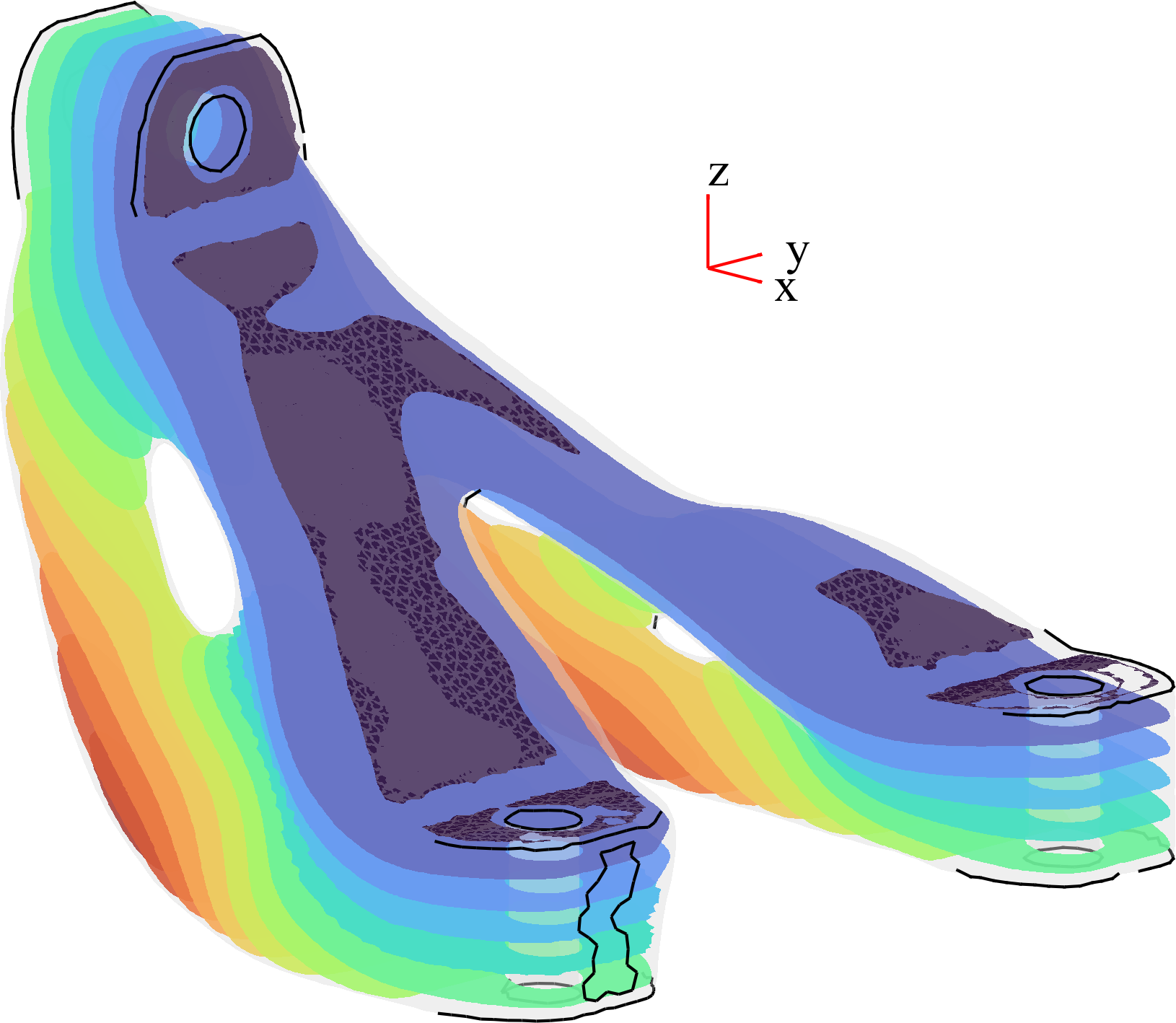}
         \caption{Offset slicing}
     \end{subfigure}
     \hfill
     \begin{subfigure}[t]{0.49\textwidth}
         \centering
         \includegraphics{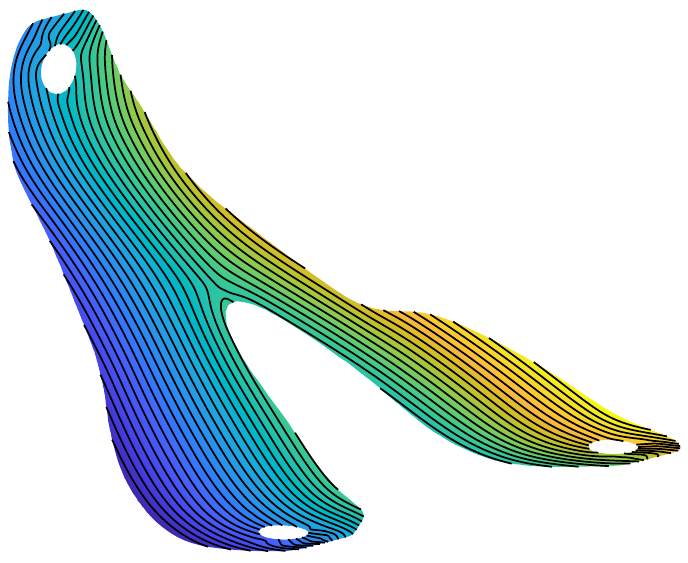}
         \caption{Isolines and scalar field on one layer of the part}
     \end{subfigure}
        \caption{Application of our method to the benchmarking piece known as \emph{topology optimization}}
        \label{fig:topopt}
\end{figure*}

Figures \ref{fig:bunny} and \ref{fig:topopt} illustrate the application of our method to two benchmarking geometries taken from the literature, which are commonly known as \emph{bunny head} and \emph{topology optimization} \cite{fang2020reinforced}.
\cmt{In Eq.\ \eqref{eq:gamma}, we defined the metric $\bar{\gamma}$ to quantify \emph{slicing} stress-alignment. Now, we are interested in \emph{trajectory} stress-alignment, for which we define a new metric.} We first compute the \cmt{trajectory} stress-alignment
\begin{equation}
    \beta_k = \|\mathbf{f} \cdot \mathbf{d}\| \,,
\end{equation}
where $k$ indicates a point on the sampled print trajectory, and $\mathbf{f}$ and $\mathbf{d}$ are the maximum principal stress and print direction at that location. We then average $\beta_k$ across all $K$ \emph{critical} points on the print trajectory to produce
\begin{equation}
    \bar{\beta} = \frac{1}{K} \sum_{\text{critical } k} \beta_k\,.
\end{equation}
The values of $\bar{\beta}$ range between $0$ (no alignment) to $1$ (perfect alignment). \cmt{In Table \ref{tab:traj_alignment}, we show the average trajectory stress-alignment $\bar{\beta}$ for the three geometries that we have demonstrated.} The results show that the trajectories produced with our method have excellent alignment with the maximum principal stress flow.

\begin{table}[htbp]
\centering
\caption{Average trajectory stress-alignment for different geometries (evaluated in the critical regions only)}
\renewcommand{\arraystretch}{1.3}
\renewcommand\cellset{\renewcommand\arraystretch{0.8}%
    \setlength\extrarowheight{0pt}}
\begin{tabular}[h]{@{}l r r r@{}}
\toprule
& Fork & Bunny Head & \makecell[c]{Topology \\ Optimization} \\
\midrule
$\bar{\beta}$ & 0.94 & 0.74 & 0.88 \\
\bottomrule
\end{tabular}
\label{tab:traj_alignment}
\end{table}

\begin{figure}[htbp]
\centerline{\includegraphics{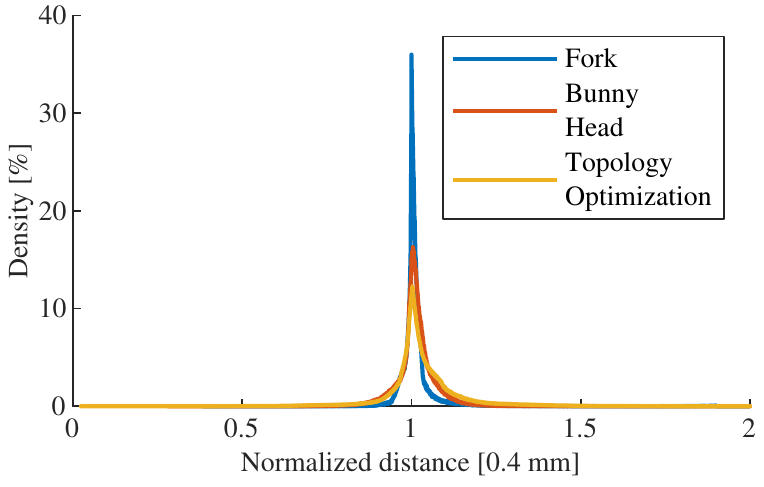}}
\caption{Distribution of distances between print trajectories for the \emph{fork}, \emph{bunny head} and \emph{topology optimization} geometries. The distance has been normalized using the nominal line spacing of \SI{0.4}{mm} that was selected when generating the trajectories.}
\label{fig:dist_hist}
\end{figure}

Finally, we show the distribution of distances between print trajectories in Fig.\ \ref{fig:dist_hist}. A tight distribution corresponds to excellent manufacturability, as the lines are very uniformly spaced. Other methods in the literature produce distributions that appear to be sums of Gaussians, with significant tails extending away from the nominal spacing (Fig.\ 13 of \cite{fang2020reinforced}). In Fig.\ \ref{fig:dist_hist}, this behavior is completely absent: the distributions are extremely tight and centered exactly on the nominal line spacing distance. Table \ref{tab:fang_comp} compares the mean and variance of the distribution of the distances produced by our method and by the method by Fang et al. \cite{fang2020reinforced}. While the means are well centered around the nominal line spacing for both methods (i.e.\ mean normalized distance $\approx 1$), the variance produced by \cite{fang2020reinforced} is worse due to the tails of the distribution. Our method, which targets manufacturability via homogeneous line spacing, exhibits instead very low distances variance. This -- paired with the homogeneous layer height produced by offset slicing -- ensures that very high-quality prints can be achieved with all materials, and particularly with sensitive ones such as LCPs.

\begin{table*}[htbp]
\centering
\caption{Mean and variance of the distribution of distances using the proposed method (see\ Fig.\ \ref{fig:dist_hist}) and the method by Fang et al \cite{fang2020reinforced}. The data from Fang et al. have been extracted from Fig.\ 13 of \cite{fang2020reinforced}. In both methods, the data have been normalized against the nominal line spacing and are thus adimensional.}
\renewcommand{\arraystretch}{1.3}
\begin{tabular}[h]{@{}l r r c r r c r r @{}}
\toprule
& \multicolumn{2}{c}{Fork} & \phantom{abc} & \multicolumn{2}{c}{Bunny Head} & \phantom{abc} & \multicolumn{2}{c}{Topology Optimization}\\
\cmidrule{2-3} \cmidrule{5-6} \cmidrule{8-9}
& Mean & Variance && Mean & Variance && Mean & Variance\\
\midrule
Proposed method & 1.01 & \num{5.0e-3} && 1.01 & \num{4.4e-3} && 1.03 & \num{13.2e-3}\\
Fang et al. \cite{fang2020reinforced} & N/A & N/A && 0.97 & \num{31.4e-3} && 0.98 & \num{31.3e-3} \\
\bottomrule
\end{tabular}
\label{tab:fang_comp}
\end{table*}

\section{Results} \label{sec:experiments}

\begin{figure}[htbp]
\centerline{\includegraphics[width=0.49\textwidth]{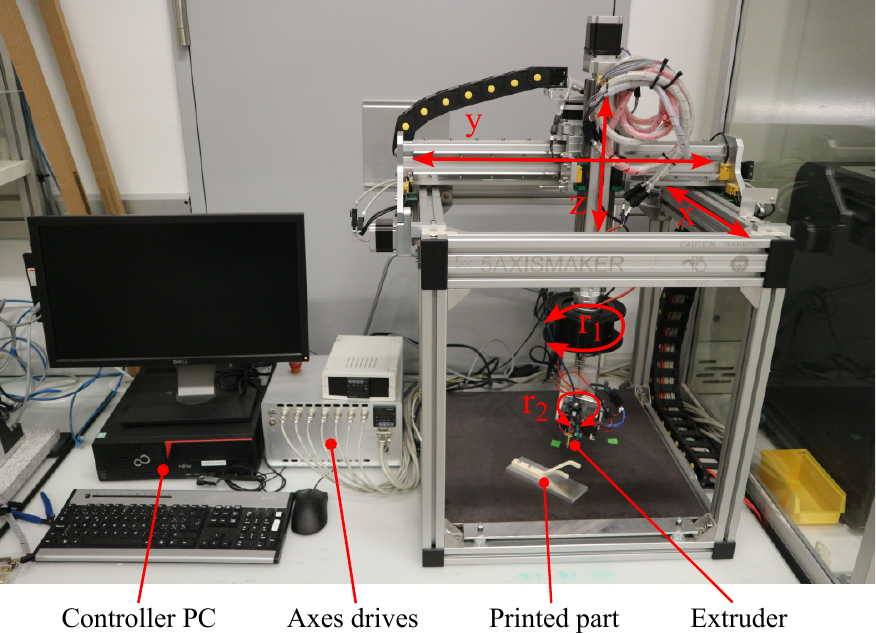}}
\caption{Custom 5-axis FFF machine with highlighted principal components}
\label{fig:printer}
\end{figure}

We have used the fork demonstrator geometry to validate experimentally the approach we have proposed and quantify the improvements achieved through stress-aligned print trajectory optimization. The parts were manufactured using a modified 5-axis computer numerical control (CNC) machine (5AxisMaker). The CNC machine was modified by replacing the milling head with a \cmtb{custom direct drive filament extruder developed by NematX AG} and was controlled using Mach3\footnote{https://www.machsupport.com/software/mach3/} on a Windows PC. Figure \ref{fig:printer} shows the custom printer and highlights its main components and axes.
\begin{figure*}[htbp]
    \begin{subfigure}[b]{0.49\textwidth}
        \centerline{\includegraphics[scale=0.9]{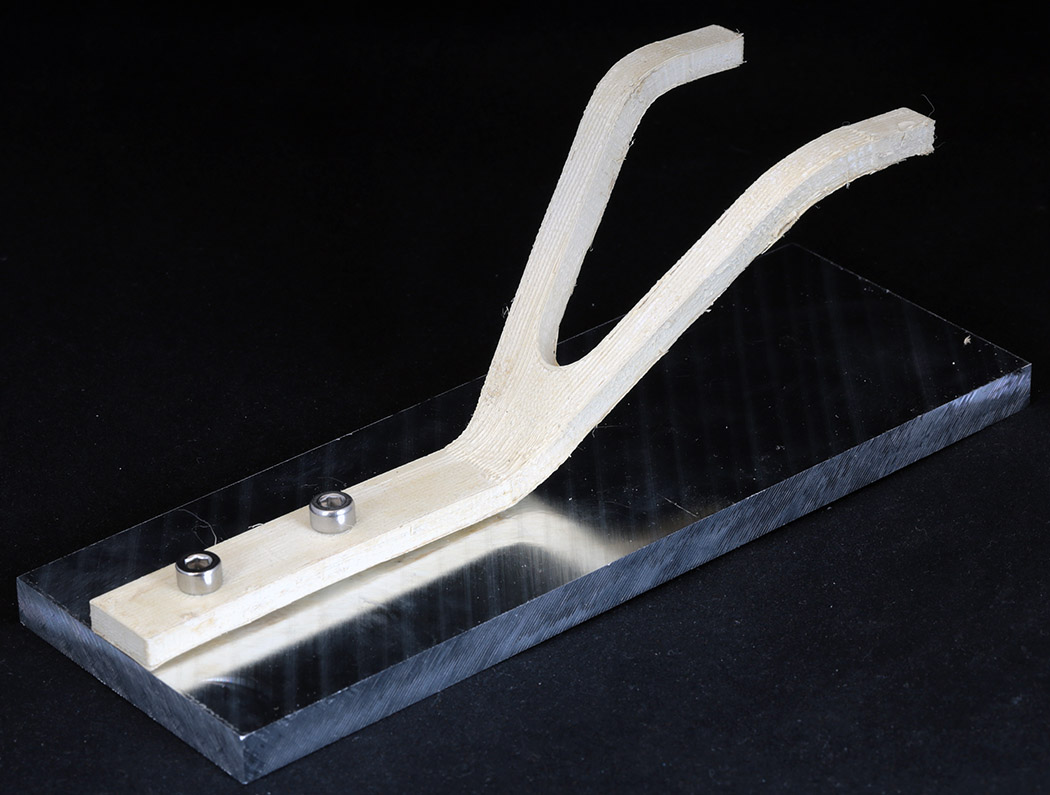}}
        \caption{Manufactured part attached to the fixture}  
    \end{subfigure}
    \hfill
    \begin{subfigure}[b]{0.49\textwidth}
        \centerline{\includegraphics[scale=0.9155]{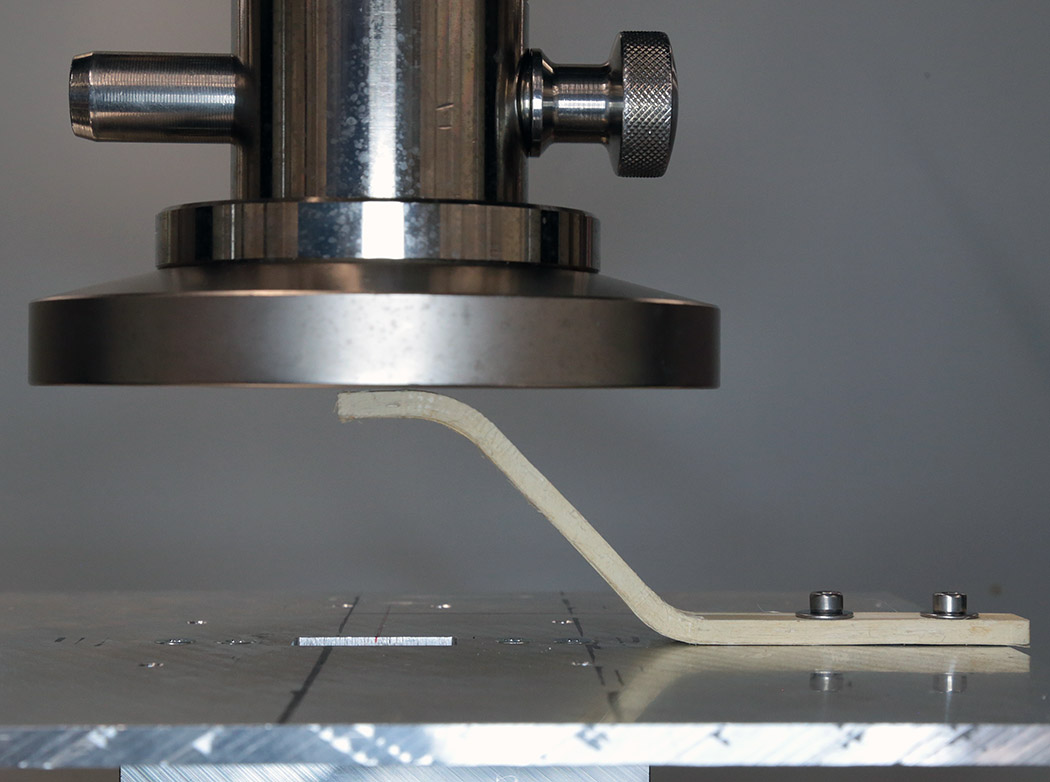}}
        \caption{Part and fixture mounted on the testing rig} 
    \end{subfigure}
    \caption{Installation of a fork manufactured using stress-aligned trajectories and LCPs on the testing rig}
    \label{fig:test_setup}   
\end{figure*}
We tested the printed samples according to the load case for which the optimization was conducted (see\ Fig.\ \ref{fig:load_case}). To do so, we used a Zwick Z020 universal testing machine with a \SI{20}{kN} capacity load cell. \cmt{As shown in Fig.\ \ref{fig:test_setup},} the fork was attached to the machine bottom sample holder with two screws passing through their holes. Then, the machine testing head was lowered at a speed of \SI{2}{\milli\metre\per\min}. The tests were stopped when the applied force dropped under a predefined threshold value, which is typically reached when the specimen is completely broken. During the test, the displacement and applied force were recorded to trace a force-displacement curve and to compute the failure force and stiffness of each printed specimen.

\subsection{Experimental Printer Impact} \label{sec:mach_impact}

\begin{figure}[htbp]
\centerline{\includegraphics[width=0.49\textwidth]{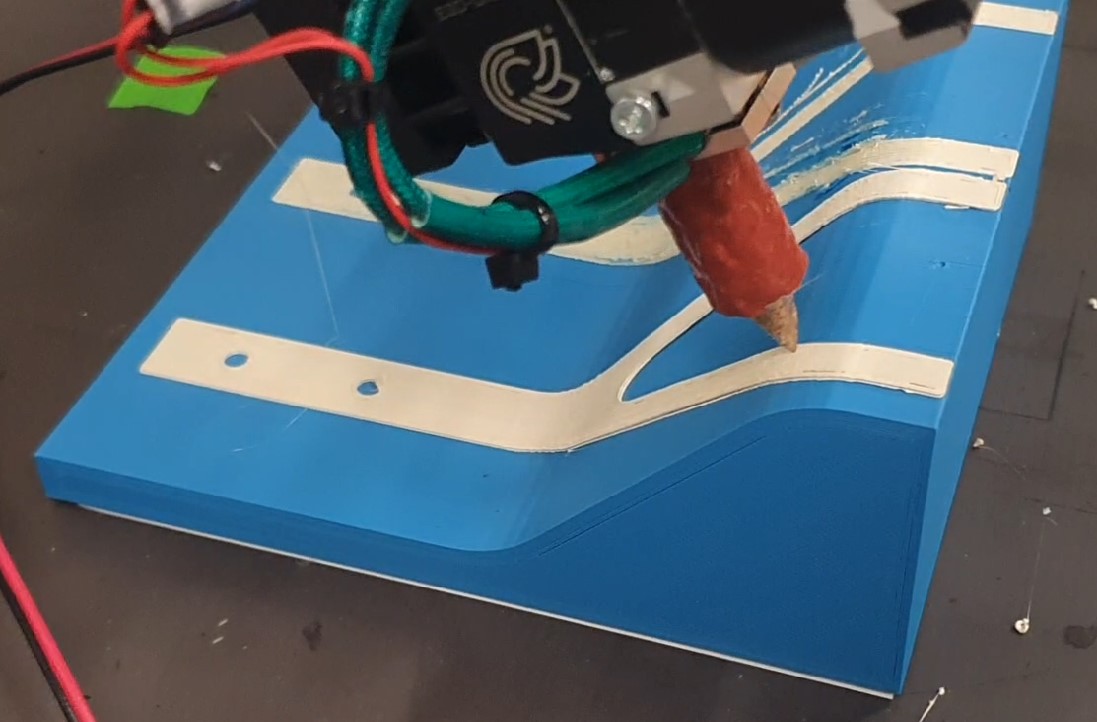}}
\caption{\cmtb{5-axis machine printing a non-planar fork on the support structure. In the 5-axis case, the shown support part was fabricated in PLA with a different machine and fixed to the 5-axis machine bed. In the 2.5D case, the Prusa i3 MK3 generated and printed the required support autonomously during the manufacturing.}}
\label{fig:support}
\end{figure}

\begin{table}[htbp]
\centering
\caption{Failure force and stiffness of the isotropic PLA samples, used to evaluate the impact of the 5-axis experimental printer on mechanical properties}
\renewcommand{\arraystretch}{1.3}
\begin{tabular}[h]{@{}l r r @{}}
\toprule
& \multicolumn{2}{c}{PLA}\\
\cmidrule{2-3}
& \makecell[c]{2.5D \\ \bf Prusa} & \makecell[c]{Non-planar \\ \bf Experimental} \\
\midrule
Failure force [\SI{}{N}] & 80.0 & 52.9 \\
Stiffness [\SI{}{N/mm}] & 6.0 & 4.6 \\
\bottomrule
\end{tabular}
\label{tab:results_PLA}
\end{table}
\cmtb{To thoroughly demonstrate the performance of our method, we plan to benchmark our optimized parts against the current state of the art in commercial LCP printing. Thus, we compare optimized non-planar parts printed by our experimental machine to conventional 2.5D parts printed by a Prusa i3 MK3s configured and used for production by NematX AG. Both machines mount the same custom extrusion assembly described in \cite{gantenbein2018three}. Print speed and temperature are shared by both machines and were also taken from \cite{gantenbein2018three}. As the extrusion assembly and configurations are identical, the differences in print quality between the two machines must be attributed to their kinematics.The \emph{self-built} experimental machine -- a modified CNC mill with a cumbersome head containing two rotational drives -- has a large print head inertia, which reduces the movement accuracy and print quality when compared to the Prusa printer \cite{avdeev2020investigation}. The 5-axis printing motion, that requires very large accelerations on all axes, further reduces the quality of non-planar prints.} Low print quality negatively affects the mechanical properties of printed components. Thus, we first evaluated the impact of our experimental 5-axis machine on part properties by printing and comparing two forks with an isotropic material (PLA). The \emph{2.5D} fork was sliced with a commercial planar slicer (Prusa Slicer) and printed on the NematX AG Prusa i3 MK3s; the \emph{non-planar} one was sliced with our proposed method and printed on the experimental 5-axis machine. Both trajectory generation methods were set to produce a layer thickness of \SI{0.1}{mm} and a line spacing (corresponding to the extrusion width) of \SI{0.4}{mm}. \cmtb{In both cases, we printed at \SI{35}{\milli\metre\per\second} and \SI{210}{\celsius} (which are the settings used in \cite{gantenbein2018three} produce PLA parts for benchmarking purposes), printed two contour lines, and used a support structure for the arms overhang \cmtb{(see Fig.\ \ref{fig:support})}.} The 2.5D slicer was set to produce an aligned rectilinear infill with 100\% density and to slice the part along the Z axis of the demonstrator bracket reference frame. \cmt{The Z axis slicing was selected for benchmarking in light of the fact that Y axis slicing requires to build a support structure on the part (between the bracket arms) and that X-axis slicing has extremely poor stress alignment (see Table \ref{tab:slice_alignment}).}
Theoretically, given the isotropic properties of PLA, the mechanical properties of the conventional and optimized samples should be similar. However, the parts manufactured with the commercial printer \cmtb{are printed more accurately}; conversely, the \cmtb{kinematics} and the experimental nature of our \cmt{self-built} printer \cmt{reduce the printing quality}. The results in Table \ref{tab:results_PLA} -- that we discuss comprehensively in Sec.\ \ref{sec:disc} -- show that the \cmtb{forks manufactured with the experimental machine during non-planar printing have worse mechanical properties.} \cmt{In the next section we will compare 2.5D and non-planar LCP samples. The machine benchmarking conducted on PLA can be used to account for the effect of the differences between the two printers on our results. This makes the difference in performance obtained via print trajectory optimization of anisotropic polymers more clearly observable.} 

\subsection{LCPs results} \label{sec:LCP_res}

\begin{figure}[htbp]
\centerline{\includegraphics{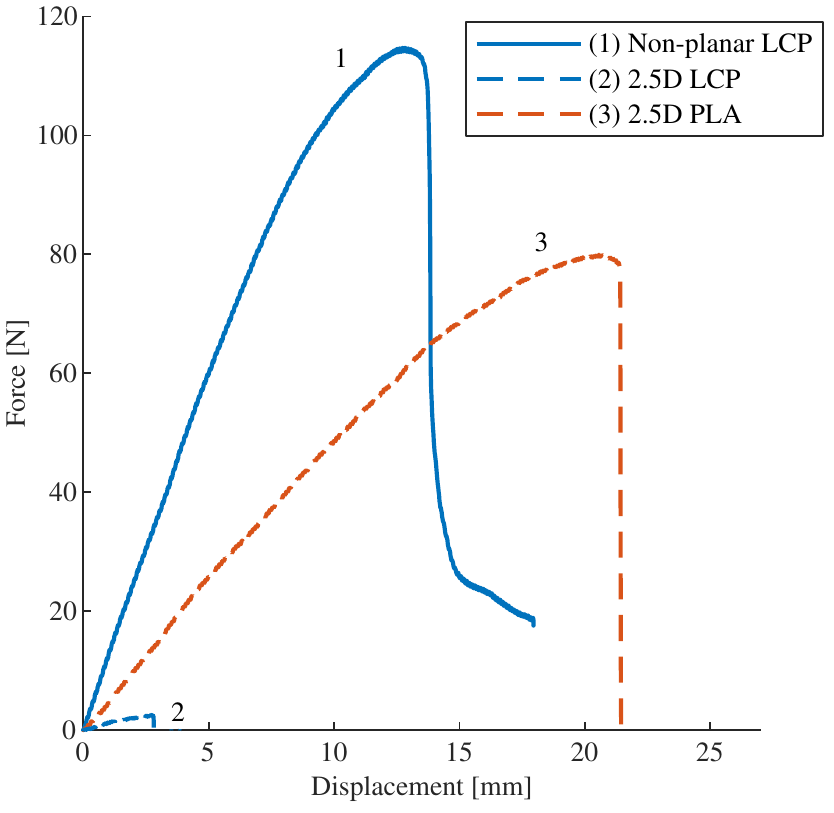}}
\caption{Force-displacement plot of the tested samples. To facilitate the comparison, the plot includes the curve of the 2.5D PLA sample produced with a conventional material, machine and slicing method.}
\label{fig:force_disp}
\end{figure}

\cmtb{Using the same machines discussed above (a Prusa i3 MK3s customized by NematX AG for commercial LCP printing and an experimental 5-axis machine, both mounting the same extrusion assembly developed in \cite{gantenbein2018three})}, we printed two samples using LCPs. The \emph{2.5D} sample was produced by planar slicing on the Prusa i3 MK3s, while the \emph{non-planar} sample was sliced with our method and printed on the 5-axis machine. \cmtb{Both samples were manufactured at \SI{295}{\celsius} at a speed of \SI{35}{\milli\metre\per\second}, and then thermally annealed for \SI{24}{\hour} at \SI{270}{\celsius}, using the state of the art print settings and annealing protocol described in \cite{gantenbein2018three}.} The results are contained in Table \ref{tab:results_LCP} and discussed in Sec.\ \ref{sec:disc}. Figure \ref{fig:force_disp} depicts the force-displacement plots of the two samples. Pictures of the LCP samples after failure are shown in Fig.\ \ref{subfig:convLCP} and \ref{subfig:optLCP}. \cmt{The machine effects discussed in the previous section affect the LCP parts as well. As the two parts have been printed on \cmtb{machines with different kinematics}, the performance of the non-planar part is worsened by the experimental 5-axis printer.}

\begin{table}[htbp]
\centering
\caption{Failure force and stiffness of the anisotropic LCP samples. The non-planar samples were sliced with the proposed optimization method and considerably outperform conventional 2.5D slicing.}
\renewcommand{\arraystretch}{1.3}
\begin{tabular}[h]{@{}l r r @{}}
\toprule
& \multicolumn{2}{c}{LCP}\\
\cmidrule{2-3} 
& \makecell[c]{2.5D \\ \bf Conventional} & \makecell[c]{Non-planar \\ \bf Optimized} \\
\midrule
Failure force [\SI{}{N}] & 2.6 & \bf 114.7 \\
Stiffness [\SI{}{N/mm}] & 2.4 & \bf 14.8 \\
\bottomrule
\end{tabular}
\label{tab:results_LCP}
\end{table}

\section{Discussion} \label{sec:disc}

The experiments demonstrate how print trajectory optimization vastly improves the mechanical properties of the parts manufactured with anisotropic materials. The failure force of the LCP samples has increased by a factor of 44, and the stiffness by a factor of 6. In light of the experimental nature of the equipment used for non-planar prints \cmtb{(whose impact on PLA prints with respect to a commercial-grade machine has been empirically quantified to -34\% failure force and -23\% stiffness in Sec.\ \ref{sec:mach_impact})}, the improvements obtained with LCP samples are even more remarkable. 

\cmt{Even when considering the performance deviations produced by repeated testing, the results would remain \cmtb{qualitatively} unchanged. For example \cmtb{-- as all parts have been produced with the materials, extrusion assembly, print settings and annealing protocol of \cite{gantenbein2018three} -- it is possible to} consider the coefficients of variation found in \cite{gantenbein2018three} (see the \SI{0}{\degree} tests of Extended Data Fig.\ 1). The machine impact (in the unfavorable scenario where the non-planar PLA sample performs \cmtb{two standard deviations} better and the 2.5D PLA sample performs \cmtb{two standard deviations} worse) would lead to a \cmtb{24\%} reduction in failure force and \cmtb{to a minor stiffness increase by a factor of 1.17}. Considering the unfavorable scenario for the LCP results (where the non-planar LCP sample performs \cmtb{two standard deviations} worse and the 2.5D LCP sample performs \cmtb{two standard deviations} better) the non-planar LCP sample would outperform the 2.5D LCP sample by a factor of \cmtb{41} in failure force and by a factor of \cmtb{5} in stiffness. In essence, this analysis suggests that the difference in performance between non-planar samples and 2.5D samples is sufficiently large to make our conclusions unaffected by process and measurement variations, even in \cmtb{a very} unfavorable scenario. \cmtb{However, given the limited amount of samples manufactured and tested for this empirical analysis, more data are necessary to provide statistically rigorous conclusions. A larger study including repeated experiments will be the subject of future work.}}

\begin{figure*}[htbp]
     \centering
     \begin{subfigure}[b]{0.49\textwidth}
         \centering
         \includegraphics[scale=0.9]{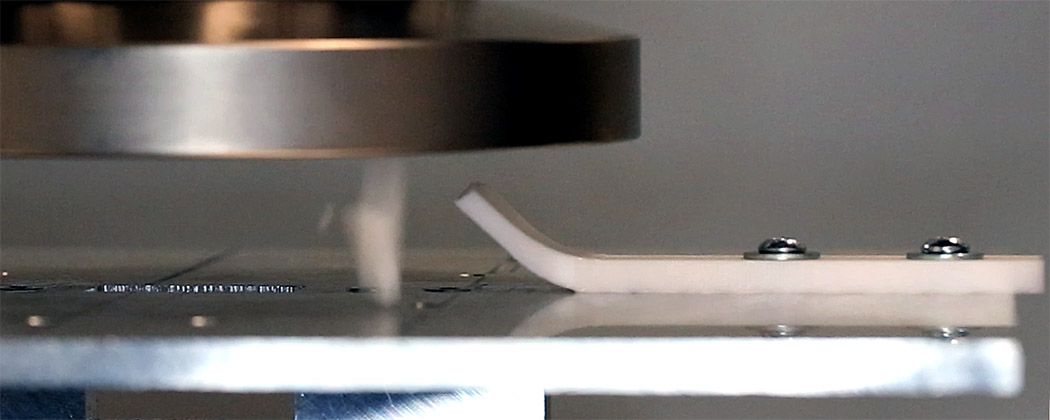}
         \caption{2.5D PLA}
         \label{subfig:convPLA}
     \end{subfigure}
     \hfill
     \begin{subfigure}[b]{0.49\textwidth}
         \centering
         \includegraphics[scale=0.9]{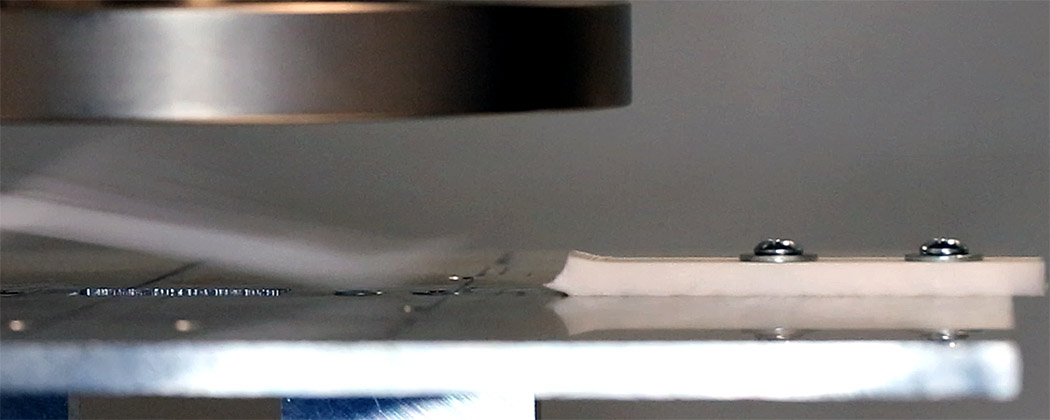}
         \caption{Non-planar PLA}
         \label{subfig:optPLA}
     \end{subfigure}
     \vskip\baselineskip
     \begin{subfigure}[b]{0.49\textwidth}
         \centering
         \includegraphics[scale=0.9]{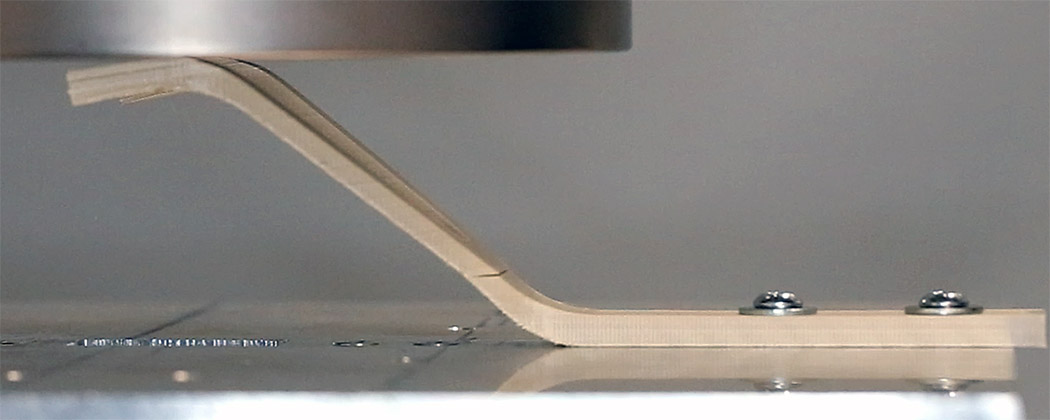}
         \caption{2.5D LCP}
         \label{subfig:convLCP}
     \end{subfigure}
     \hfill
     \begin{subfigure}[b]{0.49\textwidth}
         \centering
         \includegraphics[scale=0.288]{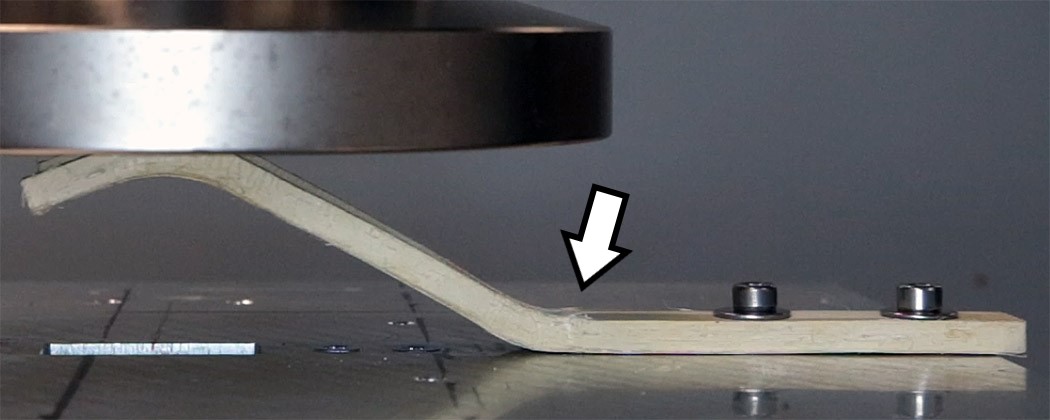}
         \caption{\cmt{Non-planar LCP}}
         \label{subfig:optLCP}
     \end{subfigure}
        \caption{Failure during testing of the four printed combinations of material and trajectory generation method}
        \label{fig:fail_samps}
\end{figure*}

Figure \ref{fig:fail_samps} offers important insights to better understand the quantitative results of the experiments. Looking first at the 2.5D prints (Fig.\ \ref{subfig:convPLA} and \ref{subfig:convLCP}), both pieces break at similar locations, in the inclined part of the bracket, where the part cross-section is minimal and horizontal layers have a very small surface. Particularly for LCPs, it is evident how the material has poor properties when the load is perpendicular to the layers' surface, producing inter-layer fractures. Thanks to its better layer adhesion properties, PLA breaks at a similar location but without a clear inter-layer failure. Moving to the case of non-planar prints, we observe how PLA (Fig.\ \ref{subfig:optPLA}) breaks at the lower curvature region. As discussed previously, PLA is quite insensitive to print orientation, and the fracture should ideally happen as in Fig.\ \ref{subfig:convPLA}. However, the lower curvature has been found to be the most complex region to print with the 5-axis equipment. The machine kinematics make high-quality material deposition in convex curvatures difficult due to large accelerations (and consequently vibrations) in the axes. It appeared upon observation that at the lower curvature the material was partially over-extruded. We believe that this defect is the cause of the mechanical properties worsening found in the non-planar PLA part, as well as for the failure location shown in Fig.\ \ref{subfig:optPLA}. Finally, Fig.\ \ref{subfig:optLCP} suggests that, despite the printing defects, the non-planar LCP fork fails at the location where FEA predicted the largest stress. Failure takes place when the deposited LCP fibers -- which are stressed in the deposition direction -- break under tension. Thanks to the optimized layers and print trajectories orientation, the part does not suffer from inter-layer or inter-trajectory fractures, and the anisotropic nature of LCPs is fully exploited to produce excellent mechanical properties.

\section{Conclusions}

In multi-axis 3D printing, it is possible to produce components using filaments that are deposited in almost any desired spatial orientation. This freedom can be exploited to improve the mechanical properties of the printed parts, particularly when using anisotropic polymers as feed stock. By aligning the material deposition direction with that of the stress flow induced by a component's load case, it is possible to dramatically improve the failure force and stiffness of the part.

In this work, we have proposed a method for non-planar slicing and print trajectory optimization that considers manufacturing constraints such as constant layer thickness and line spacing. These constraints are particularly important to obtain high-quality prints and avoid early failure of the parts. Using computations on benchmarking geometries, we showed that our method produces extremely homogeneous and thus well manufacturable trajectories when compared to existing approaches. Our framework is suitable for complex geometries and stress flows frequently encountered in common mechanical components. We have conducted experiments to manufacture load-bearing brackets using LCP material from NematX AG on a 5-axis 3D printer. The results have demonstrated that, despite the introduction of deposition defects linked with the complexity of 5-axis printing, the method fully exploits the strongly anisotropic nature of LCPs. The optimized prints achieve a $44\times$ improvement in failure strength and a $6\times$ improvement in stiffness, compared with conventional planar printing.

\section*{Funding}

This work was supported by the Swiss Innovation Agency (Innosuisse, grant \textnumero 102.617) and by the Swiss National Science Foundation under NCCR Automation (grant \textnumero 180545).

\section*{Acknowledgments}

We acknowledge the support of NematX AG that provided the LCP material and experimental setup for this work.



 \bibliographystyle{elsarticle-num} 
 \bibliography{cas-refs}





\end{document}